\documentclass[12pt]{amsart}

\usepackage{amssymb}
\usepackage{amsmath}
\usepackage{amsthm}
\usepackage{amsfonts}

\usepackage{a4wide}
\usepackage{longtable}
\usepackage{multirow}
\usepackage{xcolor}
\usepackage{hyperref}
\hypersetup{
    colorlinks=true, % make the links colored
    citecolor=blue,  % couleur des num biblio
    urlcolor=blue,   % Couleur des liens externes
}

\newtheorem{proposition}{Proposition}
\newtheorem{theorem}{Theorem}
\newtheorem{corollary}{Corollary}

\theoremstyle{definition}
\newtheorem{definition}{Definition}
\newtheorem{example}{Example}

\newtheorem{remark}{Remark}

\numberwithin{equation}{section}

\begin{document}
\title[Coclass of the second 3-class group]
{Coclass of the second 3-class group}

\author{Siham Aouissi}
\address{Number theory and information security research team\\Ecole Normale Sup\`erieure (ENS) of Moulay Ismail University (UMI)\\B.P. 3104\\Toulal\\Mekn\'es\\Morocco}
\email{s.aouissi@umi.ac.ma}
\urladdr{https://sites.google.com/view/siham-aouissi}

\author{Daniel C. Mayer}
\address{Naglergasse 53\\8010 Graz\\Austria}
\email{quantum.algebra@icloud.com}
\urladdr{http://www.algebra.at}

%\thanks{Research supported by the Austrian Science Fund (FWF): projects J0497-PHY, P26008-N25, and by the Research Executive Agency of the European Union (EUREA): Project Horizon Europe 2021--2027}
\thanks{Research of the second author supported by the Austrian Science Fund (FWF): projects J0497-PHY, P26008-N25, and by the Research Executive Agency of the European Union (EUREA)}



\date{Saturday, 23 August 2025}

\begin{abstract}
By means of parametrized presentations of finite metabelian \(3\)-groups,
it is proved that the coclass \(\mathrm{cc}(M)\)
of the second \(3\)-class group \(M=\mathrm{Gal}(\mathrm{F}_3^2(K)/K)\)
of any algebraic number field \(K\)
with elementary bicyclic \(3\)-class group \(\mathrm{Cl}_3(K)\simeq (3,3)\)
is determined unambiguously by the second largest order
\(\mathrm{ord}(\mathrm{Cl}_3(E_2))=3^{\mathrm{cc}(M)+1}\)
among the four \(3\)-class groups
of the unramified cyclic cubic extensions \(E_i\) \((i=1,\ldots,4)\) of \(K\).
Minimal discriminants of quadratic and cubic fields \(K\)
with assigned coclass \(\mathrm{cc}(M)\)
are computed from extensive databases of \(3\)-class numbers \(\mathrm{ord}(\mathrm{Cl}_3(E_i))\)
as an application.
\end{abstract}

\maketitle

\hypersetup{linkcolor=blue}

%\newpage
%--------------------------------------------------------------------------------

\section{Introduction}
\label{s:Intro}

\noindent
Let \(K/\mathbb{Q}\) be an algebraic number field.
Suppose the \(3\)-class group \(\mathrm{Cl}_3(K)=\mathrm{Syl}_3(\mathrm{Cl}(K))\),
that is the Sylow \(3\)-subgroup of the ideal class group, of \(K\) is \textit{elementary bicyclic},
\(\mathrm{Cl}_3(K)\simeq(\mathbb{Z}/3\mathbb{Z})\times(\mathbb{Z}/3\mathbb{Z})\).
This is the simplest situation where the Hilbert \(3\)-class field tower \(\mathrm{F}_3^\infty(K)\) of \(K\),
i.e., the maximal unramified pro-\(3\)-extension of \(K\),
may have a \textit{non-abelian} automorphism group \(G=\mathrm{Gal}(\mathrm{F}_3^\infty(K)/K)\).
In this case, we are interested in the \textit{coclass} \(\mathrm{cc}(M)=n-\mathrm{cl}(M)\)
of the second \(3\)-class group \(M=\mathrm{Gal}(\mathrm{F}_3^2(K)/K)\) of \(K\)
\cite{Ma2012a,Ma2013,So1933,SoTa1934},
with order \(\mathrm{ord}(M)=3^n\) and nilpotency class \(\mathrm{cl}(M)\),
that is the Galois group of the \textit{maximal metabelian} unramified \(3\)-extension \(\mathrm{F}_3^2(K)\) of \(K\),
which is the \textit{second stage} of the tower,
and is called the second Hilbert \(3\)-class field of \(K\).
According to
\cite[Cor. 3.1, p. 838]{Ma1992}
or
\cite[Thm. 3.1, p. 2214]{Ma2014b},
such a number field \(K\) has precisely \(4=\frac{3^2-1}{3-1}\) \textit{unramified cyclic cubic extensions} \(E_1,\ldots,E_4\),
and \textbf{our main result is the following convenient coclass rule}.

\begin{theorem}
\label{thm:Coclass}
The coclass \(\mathrm{cc}(M)\) of the
second \(3\)-class group \(M=\mathrm{Gal}(\mathrm{F}_3^2(K)/K)\) of an algebraic number field \(K\)
with elementary bicyclic \(3\)-class group
\(\mathrm{Cl}_3(K)\simeq(\mathbb{Z}/3\mathbb{Z})\times(\mathbb{Z}/3\mathbb{Z})\)
and four unramified cyclic cubic extensions \(E_i\) \((1\le i\le 4)\)
is given by the following relation
\begin{equation}
\label{eqn:Coclass}
\mathrm{cc}(M)=
\begin{cases}
\log_3(h_2)-1, & \text{ if } \mathrm{cl}(M)\ge 2 \text{ (non-abelian M)}, \\
\log_3(h_2), & \text{ if } \mathrm{cl}(M)=1 \text{ (abelian M)},
\end{cases}
\end{equation}
where the \(3\)-class numbers \(h_i=\#(\mathrm{Cl}_3(E_i))\)
are arranged decreasingly \(h_1\ge h_2\ge h_3\ge h_4\),
and \(\log_3\) denotes the logarithm with respect to the basis \(3\).
\end{theorem}

\noindent
The proof will be conducted in \S\
\ref{s:Main}
after preparatory sections on
\textit{abelian type and quotient invariants} in \S\
\ref{s:ATIAQI}, and
\textit{Artin pattern, IPADs, and transfer kernels and targets} in \S\
\ref{s:Artin}.
The highlights of this article are \S\
\ref{s:Arithmetic},
where the \textit{minimal absolute discriminants} or \textit{conductors}
of number fields \(K\) of degrees \(2,3,6\)
with assigned values of the coclass \(\mathrm{cc}(M)\)
are determined with immense computational effort, and \S\
\ref{s:Normal},
where \textbf{class and coclass} of \(M\) are \textbf{visualized by the body of normal lattices}.
Additional evidence of the importance of Theorems
\ref{thm:Coclass}
and
\ref{thm:Main}
is the coclass theory in \S\
\ref{s:History}.
Two crucial concepts must be explained in advance.

%--------------------------------------------------------------------------------

\subsection{The coclass of a non-abelian 3-group}
\label{ss:CoCl}

\noindent
Let \(G\) be a finite \(3\)-group with elementary bicyclic commutator quotient \(G/G^\prime\simeq (3,3)\).
The lower central series \((\gamma_i(G))_{i\ge 1}\) of \(G\) is defined by
\(\gamma_1(G)=G\) and recursively by \(\gamma_i(G)=\lbrack G,\gamma_{i-1}(G)\rbrack\), for all \(i\ge 2\).
Since \(G\) is nilpotent, there exists a unique integer \(c\ge 1\) such that
\(\gamma_c(G)>\gamma_{c+1}(G)=1\),
which is called the \textit{nilpotency class} (briefly: class) \(\mathrm{cl}(G)\) of \(G\).
If \(G\) is of order \(\mathrm{ord}(G)=3^n\), that is, of logarithmic order \(\mathrm{lo}(G)=n\),
then the \textit{coclass} of \(G\) is \(\mathrm{cc}(G)=\mathrm{lo}(G)-\mathrm{cl}(G)\).
This invariant describes the structure of the lower central series,
since the coclass is the \textit{number of bicyclic factors} in this series,
\(\mathrm{cc}(G)=\#\lbrace 1\le i\le c\mid\mathrm{ord}(\gamma_i(G)/\gamma_{i+1}(G))=9\rbrace\),
all the other factors are cyclic of order \(3\).
Another important interpretation of \(\mathrm{cc}(G)-1\) is
the \textit{number of bifurcations}
on the path starting at the abelian root \((3,3)\)
in the descendant tree \(\mathcal{G}(3)\)
of all finite \(3\)-groups and ending with \(G\):
for every increment of the coclass by one,
the path has to pass a bifurcation
between two coclass graphs \(\mathcal{G}(3,r)\) and \(\mathcal{G}(3,r+1)\).

%--------------------------------------------------------------------------------

\subsection{The bi-polarization of a non-abelian 3-group}
\label{ss:BiPol}

\noindent
Due to relations in the presentation of a
two-generated finite metabelian \(3\)-group \(M=\langle x,y\rangle\)
with \(M/M^\prime\simeq (3,3)\),
the commutator quotients \(H_i/H_i^\prime\) of the four maximal normal subgroups \(H_1,\ldots,H_4\) of \(M\)
have different orders, in general.
If \(M=\mathrm{Gal}(\mathrm{F}_3^2(K)/K)\) is the second \(3\)-class group of a number field \(K\),
then the commutator quotients \(H_i/H_i^\prime\) are
isomorphic to the \(3\)-class groups \(\mathrm{Cl}_3(E_i)\) of the four unramified cyclic cubic extensions \(E_i/K\).
Under a canonical selection of the generators \(x,y\)
\cite[Dfn. 3.2, p. 430]{Ma2014a},
the biggest order is attained by the \textit{polarization} 
\(H_1=\langle y,M^\prime\rangle\),
also called the \textit{two-step centralizer},
the second largest order by the \textit{co-polarization}
\(H_2=\langle x,M^\prime\rangle\),
and the smallest order by the remaining two subgroups
\(H_3=\langle xy,M^\prime\rangle\), \(H_4=\langle xy^{-1},M^\prime\rangle\),
which form the \textit{stabilization}.
Consequently, there is a decreasing natural arrangement
\(h_1\ge h_2\ge h_3=h_4\)
of the \(3\)-class numbers \(h_i=\#(\mathrm{Cl}_3(E_i))\).

%\newpage
%--------------------------------------------------------------------------------

\section{Abelian type invariants and abelian quotient invariants}
\label{s:ATIAQI}

%\noindent
Let \(p\ge 2\) be a prime number,
and \(A\) be a finite abelian group with order \(\mathrm{ord}(A)=p^n\) a power of \(p\) with an integer exponent \(n\ge 0\),
that is, a finite abelian \(p\)-group.

\bigskip
\noindent
\begin{proposition}
\label{prp:DirSum}
The group \(A\) is the direct sum of finite cyclic \(p\)-groups of the form
\begin{equation}
\label{eqn:DirSum}
A\simeq\left(\mathbb{Z}/p^{m_1}\mathbb{Z}\right)^{r_1}\bigoplus\ldots\bigoplus\left(\mathbb{Z}/p^{m_\sigma}\mathbb{Z}\right)^{r_\sigma}
\end{equation}
with the following uniquely determined integers:
a strictly increasing sequence of positive integer exponents \(m_1<\ldots<m_\sigma\) of length \(\sigma\ge 0\)
and associated positive integer counters of isomorphic components \(r_1,\ldots,r_\sigma\),
such that \(r_1m_1+\ldots+r_\sigma m_\sigma=n\).
The sum of the counters \(\rho=\rho_p=r_1+\ldots+r_\sigma\ge \sigma\) 
is called the \(p\)-\textbf{rank} of \(A\).
\end{proposition}

%--------------------------------------------------------------------------------

\begin{definition}
\label{dfn:ATI}
The \textit{abelian type invariants} of \(A\) are written either in \textit{power} form
\begin{equation}
\label{eqn:PowATI}
\left(\underbrace{\overbrace{p^{m_\sigma},\ldots,p^{m_\sigma}}^{r_\sigma\text{ times}},\ldots,\overbrace{p^{m_1},\ldots,p^{m_1}}^{r_1\text{ times}}}_{\rho\text{ powers}}\right)
\end{equation}
with \textbf{descending} exponents \(m_\sigma>\ldots>m_1\) or in succinct \textit{logarithmic} form 
\begin{equation}
\label{eqn:LogATI}
\mathrm{ATI}(A):=\left(m_\sigma^{r_\sigma}\ldots m_1^{r_1}\right)
\end{equation}
with counters \(r_\sigma,\ldots,r_1\) written as formal exponents indicating iteration,
and for brevity \textbf{without commas}, since the \(m_i\) will be \textbf{single digits} in this paper.
\end{definition}

\noindent
For instance, when \(p=3\), the abelian type invariants in power form
\(A\simeq (9,3,3,3,3)\), respectively \(A\simeq (81,81,27,27)\),
correspond to the logarithmic form
\(\mathrm{ATI}(A)=(21^4)\), respectively \(\mathrm{ATI}(A)=(4^23^2)\).

%--------------------------------------------------------------------------------

\noindent
Now let \(G\) be a finite \(p\)-group or
an infinite topological pro-\(p\) group with
commutator subgroup \(G^\prime\)
and finite abelianization \(G^{\mathrm{ab}}:=G/G^\prime\).

\begin{definition}
\label{dfn:AQI}
The \textit{abelian quotient invariants} of \(G\)
are the abelian type invariants of the abelianization of \(G\),
\begin{equation}
\label{eqn:AQI}
\mathrm{AQI}(G):=\mathrm{ATI}(G^{\mathrm{ab}}).
\end{equation}
\end{definition}

%--------------------------------------------------------------------------------

\begin{definition}
\label{dfn:LogOrd}
For any finite \(p\)-group \(G\) with order \(\mathrm{ord}(G)=p^n\), the \textit{logarithmic order} is
\[\mathrm{lo}(G):=n=\mathrm{cl}(G)+\mathrm{cc}(G),\ \text{the sum of class and coclass (see the introduction)}.\]
In particular, for a finite abelian \(p\)-group \(A\) with
\(\mathrm{ATI}(A)=\left(m_\sigma^{r_\sigma}\ldots m_1^{r_1}\right)\),
we have
\[\mathrm{lo}(A)=n=r_1m_1+\ldots+r_\sigma m_\sigma.\]
\end{definition}

\noindent
\textbf{Special finite abelian \(p\)-groups}:

\begin{itemize}

\item
The \textit{trivial} group, \(A=1\), has
\(\sigma=\rho=0\), and void \(\mathrm{ATI}(A)=()\).

\item
\textit{Cyclic} groups, \(A\simeq (p^{m_1})\), have
\(\sigma=\rho=r_1=1\), and \(\mathrm{ATI}(A)=(m_1)\).

\item
\textit{Homocyclic} groups, \(A\simeq (\overbrace{p^{m_1},\ldots,p^{m_1}}^{r_1\text{ times}})\), have
\(\sigma=1\), \(\rho=r_1\ge 2\), and \(\mathrm{ATI}(A)=(m_1^{r_1})\).

\item
\textbf{Nearly homocyclic} groups,
\(A\simeq\left(\underbrace{\overbrace{p^{m_2},\ldots,p^{m_2}}^{r_2\text{ times}},\overbrace{p^{m_1},\ldots,p^{m_1}}^{r_1\text{ times}}}_{\rho\text{ powers}}\right)\), have
\(\sigma=2\), \(\rho=r_1+r_2\ge 2\), and in particular \(m_2=m_1+1\),
\(\mathrm{ATI}(A)=(m_2^{r_2}m_1^{r_1})=((m_1+1)^{r_2}m_1^{r_1})\),
\(n=\mathrm{lo}(A)=r_1m_1+r_2m_2\), which becomes
\(r_1m_1+(\rho-r_1)(m_1+1)\),
and thus by Euclidean division of \(n\) by \(\rho\) with quotient \(m_1\) and non-zero remainder \(r_2\):
\[n=m_1\rho+r_2,\text{ with } 1\le r_2\le\rho-1.\]

\end{itemize}

%--------------------------------------------------------------------------------

\noindent
We shall use nearly homocyclic \(p\)-groups only in the following special case:

\begin{definition}
\label{dfn:Nearly}
For odd \(p\ge 3\), and the particular rank \(\rho=p-1\), we write
\begin{equation}
\label{eqn:Nearly} 
\mathrm{A}(p,n):=
\left(\underbrace{\overbrace{p^{m+1},\ldots,p^{m+1}}^{r_2\text{ times}},\overbrace{p^{m},\ldots,p^{m}}^{p-1-r_2\text{ times}}}_{p-1\text{ powers}}\right),
\end{equation}
also admitting zero remainder \(0\le r_2<p-1\) in the division \(0\le n=m(p-1)+r_2\).
\end{definition}

\noindent
For instance, when \(p=3\), we
include two degenerate cases \(0\le n\le 1\) with rank \(\rho\le 1\),
and we consider homocyclic groups (with even \(n\))
as special nearly homocyclic groups:

\noindent
\(\mathrm{A}(3,0)\simeq 1\,\hat{=}\,(0)\), \hfill \(\mathrm{A}(3,7)\simeq (81,27)\,\hat{=}\,(43)\), \\
\(\mathrm{A}(3,1)\simeq (3)\,\hat{=}\,(1)\), \hfill \(\mathrm{A}(3,8)\simeq (81,81)\,\hat{=}\,(4^2)\), \\
\(\mathrm{A}(3,2)\simeq (3,3)\,\hat{=}\,(1^2)\), \hfill \(\mathrm{A}(3,9)\simeq (243,81)\,\hat{=}\,(54)\),\\
\(\mathrm{A}(3,3)\simeq (9,3)\,\hat{=}\,(21)\), \hfill \(\mathrm{A}(3,10)\simeq (243,243)\,\hat{=}\,(5^2)\),\\
\(\mathrm{A}(3,4)\simeq (9,9)\,\hat{=}\,(2^2)\), \hfill \(\mathrm{A}(3,11)\simeq (729,243)\,\hat{=}\,(65)\),\\
\(\mathrm{A}(3,5)\simeq (27,9)\,\hat{=}\,(32)\), \hfill \(\mathrm{A}(3,12)\simeq (729,729)\,\hat{=}\,(6^2)\),\\
\(\mathrm{A}(3,6)\simeq (27,27)\,\hat{=}\,(3^2)\), \hfill \(\mathrm{A}(3,13)\simeq (2187,729)\,\hat{=}\,(76)\).

%\newpage
%--------------------------------------------------------------------------------

\section{The Artin pattern}
\label{s:Artin}

\noindent
Let \(p\ge 2\) be a prime number, and
\(G\) be a pro-\(p\)-group
with commutator subgroup \(G^\prime\)
and finite abelianization \(G/G^\prime\) of order \(p^v\), \(v\ge 1\).

\begin{definition}
\label{dfn:Artin}
Let
\(\mathrm{Lyr}_n(G):=\lbrace G^\prime\le H\unlhd G\mid (G:H)=p^n\rbrace\),
for \(0\le n\le v\),
be the \(v+1\) \textit{layers} of
intermediate normal subgroups between \(G\) and \(G^\prime\),
and denote by
\(T_{G,H}:\,G/G^\prime\to H/H^\prime\) the \textit{Artin transfer} homomorphism
from \(G\) to \(H\). Then we have
 
\noindent
\begin{equation}
\label{eqn:TTT}
\tau_n(G):=(\mathrm{ATI}(H/H^\prime))_{H\in\mathrm{Lyr}_n(G)}, \text{ for } 0\le n\le v,
\end{equation}
components of
the multi-layered
\textit{transfer target type} (TTT) \(\tau(G):=\lbrack\tau_0(G);\ldots;\tau_v(G)\rbrack\),

\noindent
\begin{equation}
\label{eqn:TKT}
\varkappa_n(G):=(\ker(T_{G,H}))_{H\in\mathrm{Lyr}_n(G)}, \text{ for } 0\le n\le v,
\end{equation}
components of
the multi-layered
\textit{transfer kernel type} (TKT) \(\varkappa(G):=\lbrack\varkappa_0(G);\ldots;\varkappa_v(G)\rbrack\),

\noindent
and the collection of both, \(\mathrm{AP}(G):=\left\lbrack\tau(G),\varkappa(G)\right\rbrack\),
is called the \textbf{Artin pattern} of \(G\).
The \textit{Index-\(p\) Abelianization Data} (IPAD) of \(G\),
\(\tau^{(1)}(G):=\lbrack\tau_0(G);\tau_1(G)\rbrack\),
arises by restriction to the zeroth and first layer, and
is a first order approximation of the TTT \(\tau(G)\)
\cite{BBH2017,BBH2021}.
\end{definition}

\noindent
In our applications, the prime number will always be \(p=3\)
and the abelianization of the pro-\(3\)-group \(G\)
will be elementary bicyclic \(G/G^\prime\simeq(\mathbb{Z}/3\mathbb{Z})\times(\mathbb{Z}/3\mathbb{Z})\)
of order \(3^v=3^2\) with \(v=2\).
Consequently there will be only three layers in the TTT,
the fixed bottom layer \(\tau_0(G)=\mathrm{ATI}(G/G^\prime)=(3,3)\,\hat{=}\,(1^2)\),
the crucial middle layer \(\tau_1(G)=(\mathrm{ATI}(H_1/H_1^\prime),\ldots,\mathrm{ATI}(H_4/H_4^\prime))\),
and the top layer \(\tau_2(G)=\mathrm{ATI}(G^\prime/G^{\prime\prime})\).

%\newpage
%--------------------------------------------------------------------------------

\section{The main theorem}
\label{s:Main}

\noindent
Let \(M:=G/G^{\prime\prime}\) be the metabelianization of a
pro-\(3\)-group \(G\)
with elementary bicyclic abelianization \(G/G^\prime\) of type \((3,3)\).

\begin{definition}
\label{dfn:Defect}
The \textit{lower central series} of \(M\) is defined by
\(\gamma_1(M):=M\) and
\(\gamma_j(M):=\lbrack\gamma_{j-1}(M),M\rbrack\) for \(j\ge 2\),
the \textit{nilpotency class} of \(M\),
\(c:=\mathrm{cl}(M)\),
is minimal such that
\(\gamma_c(M)>\gamma_{c+1}(M)=1\),
and the \textit{coclass} of \(M\) is \(r:=\mathrm{cc}(M)=\mathrm{lo}(M)-\mathrm{cl}(M)\).
For \(j\ge 2\),
the \(j\)th \textit{two-step centralizer} of \(M\),
\(\gamma_2(M)\le\chi_j(M)\le M\),
is maximal such that
\begin{equation}
\label{eqn:TwoStep}
\lbrack\chi_j(M),\gamma_j(M)\rbrack\le\gamma_{j+2}(M).
\end{equation}
If two isomorphism invariants are
\(s:=\min\lbrace j\ge 2\mid\chi_j(M)>\gamma_2(M)\rbrace\) with \(2\le s\le c\) and
\(e:=\min\lbrace j\ge 2\mid\gamma_{j+1}(M)/\gamma_{j+2}(M) \text{ cyclic} \rbrace\) with \(2\le e\le c\),
the \textit{defect of commutativity} \(k(M):=k\) of \(M\) is defined by \(0\le k\le 1\)
such that
\(\lbrack\chi_s(M),\gamma_e(M)\rbrack=\gamma_{c+1-k}(M)\).
\end{definition}

\begin{proposition}
\label{prp:Defect}
The \textit{defect of commutativity} \(k(M)\) of \(M\) is
given
% (phenomenologically) 
by
\begin{equation}
\label{eqn:DefCc1}
k(M)=
\begin{cases}
0 & \text{ if } M \text{ has an abelian maximal subgroup},\\
1 & \text{ if all maximal subgroups of } M \text{ are non-abelian},
\end{cases}
\end{equation}
when \(M\) is of coclass \(r=\mathrm{cc}(M)=1\), and by
\begin{equation}
\label{eqn:DefCc2}
k(M)=
\begin{cases}
0 & \text{ if the centre } \zeta_1(M) \text{ is bicyclic of type } (3,3),\\
1 & \text{ if the centre } \zeta_1(M) \text{ is cyclic of order } 3,
\end{cases}
\end{equation}
when \(M\) is of coclass \(r=\mathrm{cc}(M)\ge 2\).
\end{proposition}

%--------------------------------------------------------------------------------

\noindent
The distinction between coclass \(r=1\) and coclass \(r\ge 2\)
in Proposition \ref{prp:Defect} is necessary, since, under the given conditions,
a group of maximal class always has a cyclic centre,
and a group of non-maximal class never possesses an abelian maximal subgroup.
The first part is due to N. Blackburn \cite{Bl1958},
the second part to B. Nebelung \cite{Ne1989a,Ne1989b}.

%\newpage
%--------------------------------------------------------------------------------

\noindent
We are now ready to present our main results,
where we make use of the identifiers for finite \(3\)-groups of small order
in the SmallGroups database
\cite{BEO2002,BEO2005},
which are also implemented in the computational algebra systems
ANUPQ \cite{GNO2006},
GAP \cite{GAP2025}
and
MAGMA \cite{MAGMA2025,MAGMA6561}.

\begin{theorem}
\label{thm:Main}
(\textbf{Main Theorem on class and coclass from IPAD or TTT})
Let \(G\) be a pro-\(3\)-group
with elementary bicyclic commutator quotient \(G/G^\prime\) of type \((3,3)\)
and metabelianization \(M=G/G^{\prime\prime}\)
having nilpotency class \(c:=\mathrm{cl}(M)\) and coclass \(r:=\mathrm{cc}(M)\).

%\noindent
Then the IPAD \(\tau^{(1)}(G)=\lbrack\tau_0(G);\tau_1(G)\rbrack\) and
the TTT \(\tau(G)=\lbrack\tau_0(G);\tau_1(G);\tau_2(G)\rbrack\) of \(G\)
are given in terms of nearly homocyclic abelian \(3\)-groups by the \textbf{regular structure}
\begin{equation}
\label{eqn:IPADandTTT}
\begin{aligned}
\tau_0(G) &= \mathrm{A}(3,2);\\
\tau_1(G) &= \left(\underbrace{\overbrace{\mathrm{A}(3,c-k)}^{\text{polarization}},\overbrace{\mathrm{A}(3,r+1)}^{\text{co-polarization}}}_{\text{bi-polarization}},T_3,T_4\right);\\
\tau_2(G) &= \overbrace{\mathrm{A}(3,c-1)}^{\text{class factor}}\times\overbrace{\mathrm{A}(3,r-1)}^{\text{coclass factor}},
\end{aligned}
\end{equation}
provided that \textbf{either} \(r=1\) and \(c\ge 4\)\\
\textbf{or} \(r=2\) and \(c\ge 5\) or \(c=4\), \(k(M)=0\),\\
\textbf{or} \(r\ge 3\) and \(c\ge r+1\), except for the irregular case with \(c=r+2\), \(k(M)=1\).

\bigskip
\noindent
Here, the stabilizations (stable components) are independent of the class \(c\):
\begin{equation}
\label{eqn:Trees}
\left\lbrack T_3,T_4\right\rbrack=
\begin{cases}
\left\lbrack\mathrm{A}(3,r+1)^2\right\rbrack   & \text{ if } r=2,\ M\in\mathcal{T}^2\langle 729,54\rangle \textbf{ or } r=1,\\
\left\lbrack1^3,\mathrm{A}(3,r+1)\right\rbrack & \text{ if } r=2,\ M\in\mathcal{T}^2\langle 729,49\rangle,\\
\left\lbrack(1^3)^2\right\rbrack               & \text{ if } r=2,\ M\in\mathcal{T}^2\langle 729,40\rangle \textbf{ or } r\ge 3,
\end{cases}
\end{equation}
where \(\mathcal{T}^2R\) denotes the tree with root \(R\) and fixed coclass \(r=2\),
restricted to step size \(S=1\) (without bifurcations).
\(\mathcal{T}^2R\) is only a proper subtree of the full descendant tree \(\mathcal{T}R\).
\end{theorem}

\noindent
The missing situations with either
\(r=1\) and \(1\le c\le 3\) or
\(r=2\) and \(3\le c\le 4\)
are covered by the following corollary,
and the \textbf{irregular} case with \(c=r+2\), \(k(M)=1\)
is discussed in the next section
\ref{s:HeteroHomo}.

\begin{proof}
%The justification consists of a translation of earlier results
%into the terminology and notation of the present article.
Guided by
\cite[\S\ 2.3, pp. 476--478]{Ma2012b},
we begin with a dictionary between group theory and number theory,
since Theorem
\ref{thm:Coclass},
expressed arithmetically, is a consequence of Theorem
\ref{thm:Main},
in group theoretic formulation.
Let \(H_1,\ldots,H_4\) be the four maximal normal subgroups
of the metabelian second \(3\)-class group \(M=\mathrm{Gal}(F_3^2(K)/K)\) of \(K\),
and \(E_1,\ldots,E_4\) be the corresponding four unramified cyclic cubic extensions of \(K\), that is
\(\mathrm{Gal}(F_3^2(K)/E_i)\simeq H_i\), respectively \(E_i=\mathrm{Fix}(H_i)\).
Within the elementary bicyclic commutator quotient
\(M/M^\prime\simeq\mathrm{Cl}_3(K)\simeq (3,3)\),
which is isomorphic to the \(3\)-class group of \(K\),
the subgroups of index \(3\) are
\(H_i/M^\prime\simeq\mathrm{Norm}_{E_i/K}(\mathrm{Cl}_3(E_i))\)
isomorphic to the norm class groups of the extensions \(E_i\).
According to the Artin-Galois correpondence,
the abelian type invariants (ATI) of \(E_i\) correpond to
the abelian quotient invariants (AQI) of \(H_i\), since
\(\mathrm{Cl}_3(E_i)\simeq H_i/H_i^\prime\), for each \(1\le i\le 4\).
First, we look at the orders in \(\tau_1(G)=(\mathrm{AQI}(H_i))_{i=1}^4\),
then at their structures.

If the coclass of \(M\) is \(r=1\), then,
according to
\cite[Thm. 3.1--3.2, pp. 474--477]{Ma2012a},
\(\#(\mathrm{Cl}_3(E_1))=\#(H_1/H_1^\prime)=3^{m-k-1}=3^{c-k}\) (polarization),
since the index of nilpotency \(m\) is connected with the nilpotency class \(c\) by \(m-1=c\),
and \(\#(H_i/H_i^\prime)=3^2\), for \(i=2,3,4\) (stabilization),
except the case of abelian \(M\simeq (3,3)\),
where all \(\#(H_i/H_i^\prime)=3\), \(1\le i\le 4\).

However, if the coclass of \(M\) is \(r\ge 2\), then,
according to
\cite[Thm. 3.3--3.4, pp. 477--481]{Ma2012a},
\(\#(\mathrm{Cl}_3(E_1))=\#(H_1/H_1^\prime)=3^{m-k-1}=3^{c-k}\) (polarization),
and \(\#(\mathrm{Cl}_3(E_2))=\#(H_2/H_2^\prime)=3^{e}=3^{r+1}\) (co-polarization),
since the (bi-)cyclic factor invariant \(e\) is connected with the coclass \(r\) by \(e=r+1\),
and finally, \(\#(\mathrm{Cl}_3(E_i))=\#(H_i/H_i^\prime)=3^3\), for \(i=3,4\) (stabilization).
(For the bi-polarization, see
\cite[Dfn. 3.2, p. 430]{Ma2013}.)

Now we come to the \textit{structure} of the abelian \(3\)-groups
\(\mathrm{Cl}_3(E_i)\simeq H_i/H_i^\prime\), for \(1\le i\le 4\).

If the coclass of \(M\) is \(r=1\), then,
according to
\cite[Thm. 3.1, p. 421]{Ma2014a},
\(\mathrm{Cl}_3(E_1)\simeq\mathrm{A}(3,c-k)\), when \(c=m-1\ge 5\), and
\(\mathrm{Cl}_3(E_i)\simeq\mathrm{A}(3,2)\), when \(i=2,3,4\) and \(c\ge 3\),
i.e., for sufficiently large class,
which has to be supplemented by Corollary
\ref{cor:Exceptions},
for \(c\le 4\).

However, if the coclass of \(M\) is \(r\ge 2\), then,
according to
\cite[Thm. 3.2, p. 424]{Ma2014a},
\(\mathrm{Cl}_3(E_1)\simeq\mathrm{A}(3,c-k)\), for \(c=m-1\ge 5\), and
\(\mathrm{Cl}_3(E_2)\simeq\mathrm{A}(3,r+1)\), for \(r=e-1\ge 3\),
again only for sufficiently large class and coclass. The decision
\(\mathrm{Cl}_3(E_i)\simeq\mathrm{A}(3,3)\) or \(\simeq (3,3,3)\), when \(i=3,4\),
is settled according to
\cite[Thm. 4.4, p. 440, Tbl. 4.7, p. 441]{Ma2014a},
for coclass \(r=2\) by three coclass trees.
According to
\cite[Thm. 4.5, pp. 444--445]{Ma2014a},
the claims are true for large coclass \(r\ge 3\).
\(\tau_0(G)\) is trivial and
\(\tau_2(G)\) will be supplied by Theorem
\ref{thm:Commutator}.
Also compare with
\cite[Main Theorem 2.2, p. 16]{Ma29JA}, and
\cite[Thm. 3.1--3.2, pp. 290--291]{Ma2015b}.
\end{proof}

%--------------------------------------------------------------------------------

\begin{corollary}
\label{cor:Exceptions}
(\textbf{Exceptions from the Main Theorem})
The remaining groups, which are not covered by Theorem \ref{thm:Main},
share the common bottom layer \(\tau_0(G)=\mathrm{A}(3,2)\),
and possess \textbf{partially irregular}
middle layer \(\tau_1(G)\) and top layer \(\tau_2(G)\)
as given in Table \ref{tbl:Exceptions}.
\end{corollary}

\renewcommand{\arraystretch}{1.1}

\begin{table}[ht]
\caption{Regular and irregular IPAD and TTT}
\label{tbl:Exceptions}
\begin{center}

{\small

\begin{tabular}{|lccc||cc|cc|}
\hline
  \(G\)                 & \(c\) & \(r\) & \(k\) & \(\tau_1(G)\)                                                                     & Remark    & \(\tau_2(G)\)                            & Remark \\
\hline
 \(\langle 9,2\rangle\)             & 1 & 1 & 0 & \(\lbrack\mathrm{A}(3,1),\mathrm{A}(3,1),\mathrm{A}(3,1),\mathrm{A}(3,1)\rbrack\) & irregular & \(\mathrm{A}(3,0)\times\mathrm{A}(3,0)\) & regular \\
 \(\langle 27,3\rangle\)            & 2 & 1 & 0 & \(\lbrack\mathrm{A}(3,2),\mathrm{A}(3,2),\mathrm{A}(3,2),\mathrm{A}(3,2)\rbrack\) & regular   & \(\mathrm{A}(3,1)\times\mathrm{A}(3,0)\) & regular \\
 \(\langle 27,4\rangle\)            & 2 & 1 & 0 & \(\lbrack\mathrm{A}(3,2),(2),(2),(2)\rbrack\)                                     & irregular & \(\mathrm{A}(3,1)\times\mathrm{A}(3,0)\) & regular \\
 \(\langle 81,7\rangle\)            & 3 & 1 & 0 & \(\lbrack 1^3,\mathrm{A}(3,2),\mathrm{A}(3,2),\mathrm{A}(3,2)\rbrack\)            & irregular & \(\mathrm{A}(3,2)\times\mathrm{A}(3,0)\) & regular \\
 \(\langle 81,8\ldots 10\rangle\)   & 3 & 1 & 0 & \(\lbrack\mathrm{A}(3,3),\mathrm{A}(3,2),\mathrm{A}(3,2),\mathrm{A}(3,2)\rbrack\) & regular   & \(\mathrm{A}(3,2)\times\mathrm{A}(3,0)\) & regular \\
 \(\langle 243,28\ldots 30\rangle\) & 4 & 1 & 1 & \(\lbrack\mathrm{A}(3,3),\mathrm{A}(3,2),\mathrm{A}(3,2),\mathrm{A}(3,2)\rbrack\) & regular   & \(\mathrm{A}(3,3)\times\mathrm{A}(3,0)\) & regular \\ 
\hline
 \(\langle 243,3\rangle\)           & 3 & 2 & 0 & \(\lbrack\mathrm{A}(3,3),\mathrm{A}(3,3),(1^3),(1^3)\rbrack\)                     & regular   & \(\mathrm{A}(3,2)\times\mathrm{A}(3,1)\) & regular \\
 \(\langle 729,34\ldots 36\rangle\) & 4 & 2 & 1 & \(\lbrack\mathrm{A}(3,3),\mathrm{A}(3,3),(1^3),(1^3)\rbrack\)                     & regular   & \(\mathrm{A}(3,2)\times\mathrm{A}(3,2)\) & irregular \\
 \(\langle 729,37\ldots 39\rangle\) & 4 & 2 & 1 & \(\lbrack\mathrm{A}(3,3),\mathrm{A}(3,3),(1^3),(1^3)\rbrack\)                     & regular   & \(\mathrm{A}(3,3)\times\mathrm{A}(3,1)\) & regular \\
 \(\langle 243,4\rangle\)           & 3 & 2 & 0 & \(\lbrack (1^3),(1^3),\mathrm{A}(3,3),(1^3)\rbrack\)                              & irregular & \(\mathrm{A}(3,2)\times\mathrm{A}(3,1)\) & regular \\
 \(\langle 729,44\ldots 47\rangle\) & 4 & 2 & 1 & \(\lbrack (1^3),(1^3),\mathrm{A}(3,3),(1^3)\rbrack\)                              & irregular & \(\mathrm{A}(3,3)\times\mathrm{A}(3,1)\) & regular \\
 \(\langle 243,5\rangle\)           & 3 & 2 & 0 & \(\lbrack\mathrm{A}(3,3),\mathrm{A}(3,3),(1^3),\mathrm{A}(3,3)\rbrack\)           & irregular & \(\mathrm{A}(3,2)\times\mathrm{A}(3,1)\) & regular \\
 \(\langle 243,6\rangle\)           & 3 & 2 & 0 & \(\lbrack\mathrm{A}(3,3),\mathrm{A}(3,3),(1^3),\mathrm{A}(3,3)\rbrack\)           & regular   & \(\mathrm{A}(3,2)\times\mathrm{A}(3,1)\) & regular \\
 \(\langle 243,7\rangle\)           & 3 & 2 & 0 & \(\lbrack (1^3),\mathrm{A}(3,3),(1^3),\mathrm{A}(3,3)\rbrack\)                    & irregular & \(\mathrm{A}(3,2)\times\mathrm{A}(3,1)\) & regular \\
 \(\langle 243,8\rangle\)           & 3 & 2 & 0 & \(\lbrack\mathrm{A}(3,3),\mathrm{A}(3,3),\mathrm{A}(3,3),\mathrm{A}(3,3)\rbrack\) & regular   & \(\mathrm{A}(3,2)\times\mathrm{A}(3,1)\) & regular \\
 \(\langle 243,9\rangle\)           & 3 & 2 & 0 & \(\lbrack\mathrm{A}(3,3),\mathrm{A}(3,3),\mathrm{A}(3,3),\mathrm{A}(3,3)\rbrack\) & irregular & \(\mathrm{A}(3,2)\times\mathrm{A}(3,1)\) & regular \\
 \(\langle 729,56\ldots 57\rangle\) & 4 & 2 & 1 & \(\lbrack\mathrm{A}(3,3),\mathrm{A}(3,3),\mathrm{A}(3,3),\mathrm{A}(3,3)\rbrack\) & irregular & \(\mathrm{A}(3,2)\times\mathrm{A}(3,2)\) & irregular \\
\hline
\end{tabular}

}

\end{center}
\end{table}

\begin{proof}
We take information and identifiers from the SmallGroups library \cite{BEO2005},
and we begin with coclass \(r=1\) and class \(1\le c\le 3\).
The behavior of the abelian root \(\langle 9,2\rangle\)
and of the groups \(\langle 27,4\rangle\) and \(\langle 81,7\rangle\) is irregular,
whereas \(\langle 27,3\rangle\) and \(\langle 81,8\ldots 10\rangle\) fit into 
the rules of Formula \eqref{eqn:IPADandTTT}.
By CT, we abbreviate the \textit{capitulation type} (TKT component \(\varkappa=\varkappa_1\)).

Independently of \cite{BEO2005}, if the coclass of \(M\) is \(r=1\), then,
for small class \(c\le 4\),
\cite[Thm. 4.1, p. 427, Cor. 4.1.1, p. 428, Tbl. 4.1, p. 429]{Ma2014a}
clarifies the irregular structures
\(\mathrm{Cl}_3(E_i)\simeq (9)\), \(i=2,3,4\), for \(\langle 27,4\rangle\simeq G_0^3(0,1)\) with \(c=2\),
and \(\mathrm{Cl}_3(E_1)\simeq (3,3,3)\) for \(\langle 81,7\rangle\simeq G_0^4(1,0)\simeq\mathrm{Syl}_3A_9\) with \(c=3\),
whereas \(\langle 27,3\rangle\simeq G_0^3(0,0)\) with \(c=2\),
\(\langle 81,8\rangle\simeq G_0^4(-1,0)\),
\(\langle 81,9\rangle\simeq G_0^4(0,0)\) and
\(\langle 81,10\rangle\simeq G_0^4(0,1)\) with \(c=3\), and
\(\langle 243,28\rangle\simeq G_1^5(0,-1)\),
\(\langle 243,29\rangle\simeq G_1^5(0,1)\), and
\(\langle 243,30\rangle\simeq G_1^5(0,0)\) with \(c=4\), \(k(M)=1\),
behave regularly.

Now we come to coclass \(r=2\) and class \(3\le c\le 4\).
All sporadic groups outside of coclass trees are necessarily irregular,
since Formula \eqref{eqn:Trees} cannot be applied to them.
However, if this formula is extended to parents, then some exceptional groups can be viewed as regular, since
\(\langle 243,3\rangle\) with CT b.10, \(\varkappa\sim (0043)\), is the parent of the periodic root of the coclass tree \(\mathcal{T}^2\langle 729,40\rangle\);
\(\langle 243,6\rangle\) with CT c.18, \(\varkappa\sim (0313)\), is parent of the root of \(\mathcal{T}^2\langle 729,49\rangle\); and
\(\langle 243,8\rangle\) with CT c.21, \(\varkappa\sim (0231)\), is parent of the root of \(\mathcal{T}^2\langle 729,54\rangle\).
The top layer \(\tau_2(G)\) is irregular (homocyclic) only for \(\langle 729,34\ldots 36\rangle\) and \(\langle 729,56\ldots 57\rangle\).
Irregular structures are clarified by
\cite[Thm. 4.2, p. 434, Tbl. 4.3, p. 434]{Ma2014a}
for the top vertices of the coclass graph \(\mathcal{G}(3,2)\)
with order \(3^5=243\) and class \(c=3\), i.e. the \textit{isoclinism family} \(\Phi_6\),
namely \(\langle 243,4\rangle\) with CT H.4, \(\varkappa\sim (4443)\),
\(\langle 243,5\rangle\) with CT D.10, \(\varkappa\sim (2241)\),
\(\langle 243,7\rangle\) with CT D.5, \(\varkappa\sim (4224)\), and
\(\langle 243,9\rangle\) with CT G.19, \(\varkappa\sim (2143)\).
Further by
\cite[Thm. 4.3, p. 438, Tbl. 4.5, p. 438]{Ma2014a},
we get information
for the isoclinism families \(\Phi_t\), \(40\le t\le 43\), as follows.
CT b.10, \(\varkappa\sim (0043)\):
\(\langle 729,i\rangle\) with irregular \(\tau_2=(1^4)\) for \(i\in\lbrace 34,35,36\rbrace\), family \(\Phi_{40}\);
with regular \(\tau_2=(21^2)\) for \(i\in\lbrace 37,38,39\rbrace\), family \(\Phi_{41}\).
CT H.4, \(\varkappa\sim (4443)\):
\(\langle 729,i\rangle\) with regular \(\tau_2=(21^2)\) for \(i\in\lbrace 44,45,46,47\rbrace\), family \(\Phi_{42}\).
CT G.19, \(\varkappa\sim (2143)\):
\(\langle 729,i\rangle\) with irregular \(\tau_2=(1^4)\) for \(i\in\lbrace 56,57\rbrace\), family \(\Phi_{43}\).
\end{proof}

%--------------------------------------------------------------------------------

\begin{remark}
\label{rmk:Exceptions}
Illuminating tree diagrams of \(3\)-groups \(G\)
which occur in the proof of Corollary
\ref{cor:Exceptions}
are drawn
for coclass \(\mathrm{cc}(G)=1\) in
\cite[Fig. 1, p. 44]{Ma2018a}
and
\cite[Fig. 1--2, pp. 142--143]{Ma2017}
for coclass \(\mathrm{cc}(G)=2\) in
\cite[Fig. 3--5, pp. 151--153]{Ma2017}.
\end{remark}

%--------------------------------------------------------------------------------

\begin{proof}
(\textbf{Proof} of Theorem \ref{thm:Coclass})
In contrast to the proof of Theorem
\ref{thm:Main},
where the precise structure of the four \(3\)-class groups \(\mathrm{Cl}_3(E_i)\) is specified,
we now only need \(3\)-class numbers, i.e., orders of \(3\)-class groups.
Let the \(3\)-class numbers of
the four cyclic cubic extensions \(E_i\) of the base field \(K\) be
\(h_i=\mathrm{ord}(\mathrm{Cl}_3(E_i))\), for \(1\le i\le 4\),
and let the nilpotency class, coclass and the defect of commutativity
of the second \(3\)-class group \(M=\mathrm{Gal}(\mathrm{F}_3^2(K)/K)\)
be \(c=\mathrm{cl}(M)\), \(r=\mathrm{cc}(M)\), and \(k\).
Then, due to the distinguished role of the polarization,
the \(3\)-class numbers are automatically arranged in decreasing order,
\(h_1=3^{c-k}\ge h_2=h_3=h_4=3^2\), for coclass \(r=1=\log_3(h_2)-1\)
\cite[Thm. 3.1, pp. 474--475, Thm. 3.2, p. 477]{Ma2015a},
since \(k=0\) for \(2\le c\le 3\), and \(0\le k\le 1\) for \(c\ge 4\),
whence \(c-k\ge 2\).
The single exception is the abelian group
\(M=(\mathbb{Z}/3\mathbb{Z})\times(\mathbb{Z}/3\mathbb{Z})\)
with \(c=1\),
where \(h_1=h_2=h_3=h_4=3\), and thus the coclass is \(r=1=\log_3(h_2)\).
(It is illuminating to compare with
\cite[Thm. 3.2, p. 416]{Ma2013},
and to observe the supplement
\cite[Thm 3.5, p. 420]{Ma2013}.)
On the other hand,
due to the distinguished role of the bi-polarization,
\(h_1=3^{c-k}\ge h_2=3^{r+1}\ge h_3=h_4=3^3\), for coclass \(r\ge 2\)
\cite[Thm. 3.3, pp. 477--478, Thm. 3.4, p. 481]{Ma2015a}.
For \(r\ge 2\), we also have \(0\le k\le 1\) for \(c\ge r+2\),
and \(k=0\) for \(c=r+1\),
so \(r+1\ge 3\) and \(c-k\ge r+1\),
because a path from the abelian root \(\langle 9,2\rangle\)
to a group of minimal order with assigned coclass \(r\)
must consist of a single step of size \(S=1\)
to the group \(\langle 27,3\rangle\) of class \(2\)
and then \(r\) steps of size \(S=2\),
which gives at least class  \(1+r\), that is, \(c\ge r+1\),
but also \(c-k\ge r+1\), since, for interface groups, \(c=r+1\) implies \(k=0\).
It is illuminating to compare with
\cite[Thm. 3.12, pp. 435--436]{Ma2013},
and to observe the supplement
\cite[Thm 3.13, pp. 436--437]{Ma2013}.
\end{proof}

%\newpage
%--------------------------------------------------------------------------------

\section{Heterocyclic and homocyclic commutator groups}
\label{s:HeteroHomo}

\noindent
According to
\cite{Ne1989a}
and
\cite{Ma2018b},
the commutator subgroup \(G^\prime\)
of a metabelian \(3\)-group \(G\)
with elementary bicyclic abelianization \(G/G^\prime\)
and coclass \(\mathrm{cc}(G)\ge 2\),
is usually a heterocyclic direct product of two nearly homocyclic abelian \(3\)-groups,
as given in the \textbf{regular} scenario of Theorem
\ref{thm:Commutator}.
For certain immediate descendants of interface groups, however,
the structure of \(G^\prime\) may degenerate to
the direct product of two identical nearly homocyclic abelian \(3\)-groups of even logarithmic order,
whence \(G^\prime\) is in fact homocyclic.
This situation is called \textbf{irregular}.

The criterion for a homocyclic commutator subgroup \(G^\prime\)
cannot be expressed in terms of class and coclass alone.
Therefore we must use a normalized and parametrized power-commutator-presentation of \(G\).
Generators \(x,y\) of \(G=\langle x,y\rangle\) are selected such that
the elementary bicyclic third factor
\(\gamma_3(G)=\langle x^3,y^3,\gamma_4(G)\rangle\)
of the lower central series of \(G\)
is generated by the powers \(x^3,y^3\),
whereas \((xy)^3,(xy^{-1})^3\) are contained deeper
in the second center \(\zeta_2(G)\)
of the upper central series of \(G\).
The main commutator is \(s_2=t_2=\lbrack y,x\rbrack\in\gamma_2(G)\).
Moreover, \(y\in\chi_s(G)\setminus\gamma_2(G)\)
is required to lie in the first non-trivial two-step-centralizer,
whereas \(x\in G\setminus\chi_s(G)\) is chosen outside,
provided that \(G\) is not an interface group with \(c=r+1\).
Under these assumptions,
four series of commutators are defined by
\(s_j=\lbrack s_{j-1},x\rbrack\) for \(3\le j\le c+1\),
\(t_j=\lbrack t_{j-1},y\rbrack\) for \(3\le j\le r+3\),
\(\sigma_3=y^3\), \(\sigma_j=\lbrack\sigma_{j-1},x\rbrack\) for \(4\le j\le c+1\), and
\(\tau_3=x^3\), \(\tau_j=\lbrack\tau_{j-1},y\rbrack\) for \(4\le j\le r+3\).
Then \(s_j,t_j,\sigma_j,\tau_j\in\gamma_j(G)\), and
the nilpotency of \(G\) is expressed by \(s_{c+1}=\sigma_{c+1}=t_{r+3}=\tau_{r+3}=1\).
With exponents \(-1\le\alpha,\beta,\gamma,\delta,\varrho\le 1\) as parameters,
the following relations are satisfied
\(s_2^3=\sigma_4\sigma_c^{-\varrho\beta}\tau_4^{-1}\),
\(s_3\sigma_3\sigma_4=\sigma_{c-1}^{\varrho\beta}\sigma_c^\gamma\tau_{r+1}^\delta\),
\(t_3^{-1}\tau_3\tau_4=\sigma_{c-1}^{\varrho\delta}\sigma_c^\alpha\tau_{r+1}^\beta\),
\(\tau_{r+2}=\sigma_c^{-\varrho}\). See
\cite{Ne1989a}
and
\cite[\S\ 8, pp. 456--461]{Ma2014a}.

\begin{theorem}
\label{thm:Commutator}
The structure of the abelian commutator subgroup \(M^\prime\)
of the metabelianization \(M=G/G^{\prime\prime}\),
with order \(\mathrm{ord}(M)=3^n\), class \(c=\mathrm{cl}(M)\), and coclass \(r=\mathrm{cc}(M)\),
of a pro-\(3\) group \(G\),
with elementary bicyclic abelianization \(G/G^\prime\),
is given by the following direct product of nearly homocyclic abelian \(3\)-groups:
\begin{equation}
\label{eqn:Commutator}
M^\prime\simeq
\begin{cases}
\mathrm{A}(3,c-2)\times\mathrm{A}(3,c-2) & \text{ in the irregular \textbf{homocyclic} case,}\\
\mathrm{A}(3,c-1)\times\mathrm{A}(3,r-1)  & \text{ in the regular \textbf{heterocyclic} case.}
\end{cases}
\end{equation}
The irregular case is characterized by
even logarithmic order \(n=2c-2\equiv 0\,(\mathrm{mod}\,2)\),
even coclass \(r=c-2\equiv 0\,(\mathrm{mod}\,2)\),
even class \(c\equiv 0\,(\mathrm{mod}\,2)\),
and relational parameters, \(\varrho,\beta\),
\begin{equation}
\label{eqn:Rho}
\begin{cases}
0\ne\varrho=\beta-1 & \text{ for coclass } r=2,\\
\varrho=-1          & \text{ for coclass } r\ge 4.
\end{cases}
\end{equation}
\end{theorem}

\noindent
This is
\cite[Thm. 8.8, p. 461]{Ma2014a}.
It follows from the mandatory inequality \(c\ge r+1\) for non-abelian \(3\)-groups
that \(\mathrm{A}(3,c-1)\times\mathrm{A}(3,r-1)\) is really \textbf{hetero}cyclic,
and can never be homocyclic.
In the next section \S\ \ref{s:Arithmetic},
we use the translation from group theory to number theory in
\cite[\S\ 2.3, pp. 476--478]{Ma2012b}.

%\newpage
%--------------------------------------------------------------------------------

\section{Arithmetical realization of groups with assigned coclass}
\label{s:Arithmetic}

\noindent
In order to obtain arithmetical realizations
of a finite metabelian \(3\)-group \(M\) with assigned coclass \(\mathrm{cc}(M)=r\)
as the second \(3\)-class group \(M\simeq\mathrm{Gal}(\mathrm{F}_3^2(K)/K)\)
of an algebraic number field \(K\),
we first selected \(K=\mathbb{Q}(\sqrt{d})\) to be a quadratic number field
with discriminant \(d\) and elementary bicyclic \(3\)-class group.
We determined all IPADs \(\lbrack\tau_0(K);\tau_1(K)\rbrack\) in the range \(-10^6<d<10^7\)
\cite{Ma2012a,Ma2014a},
using the class field theoretic routines by Fieker \cite{Fi2001},
which are implemented in the computational algebra system Magma
\cite{MAGMA2025}.
It suffices to restrict to the middle layer
\(\tau_1(K)=(\mathrm{ATI}(\mathrm{Cl}_3(E_i)))_{i=1}^4\),
simply denoted by \(\tau(K)\),
since the bottom component \(\tau_0(K)=(11)\) remains fixed
throughout this article.
See Table
\ref{tbl:Imag} and \ref{tbl:Real}.
Only the \textit{biggest entry} in Table
\ref{tbl:Imag}, respectively \ref{tbl:Real},
is based on other, most extensive, lists of IPADs computed by M. R. Bush for
\cite{BBH2017}, respectively \cite{BBH2021},
with \(-10^8<d<10^9\).

%\newpage
%--------------------------------------------------------------------------------

%\section{Main Theorem}
%\label{s:MainX}

%\noindent

%\begin{theorem}
%\label{thm:MainX}
%\begin{equation}
%\label{eqn:MainX}
%\begin{aligned}
%\end{aligned}
%\end{equation}
%\end{theorem}

%\begin{proof}
%\end{proof}

%\newpage
%--------------------------------------------------------------------------------

\subsection{Imaginary quadratic fields}
\label{ss:Imag}

\noindent
The coclass \(\mathrm{cc}(M)\) of the
second \(3\)-class group \(M\simeq\mathrm{Gal}(\mathrm{F}_3^2(K)/K)\)
of an imaginary quadratic number field \(K\) with negative discriminant \(d<0\)
must be an \textit{even} positive integer
\cite{Ma2012a,Ma2013,Ma2014a}.
Table \ref{tbl:Imag} shows the minimal absolute discriminants \(\lvert d\rvert\)
of imaginary quadratic number fields \(K\) with assigned even coclass \(\mathrm{cc}(M)\),
their prime factorization, the capitulation type (CT) \(\varkappa(K)\),
and, most important, the logarithmic abelian type invariants, collected in the IPAD \(\tau(K)\),
where the crucial \textbf{co-polarization} \(E_2\) with
\(\mathrm{lo}(\mathrm{Cl}_3(E_2))=\mathrm{cc}(M)+1\)
is emphasized by \textbf{boldface} font.
For instance, in the last row, using the identification \(\mathbf{(54)}\,\hat{=}\,\lbrack 3^5,3^4\rbrack=\lbrack 243,81\rbrack\),
we have
\[
\text{cc}(M)+1=\mathbf{8}+1=9=5+4=\mathrm{lo}(\mathbf{(54)})=\mathrm{lo}(\lbrack 243,81\rbrack)=\mathrm{lo}(\mathrm{Cl}_3(E_2)).
\]

\renewcommand{\arraystretch}{1.1}

\begin{table}[ht]
\caption{Imaginary quadratic number fields \(K\), OEIS A380102 \cite{OEIS2025}}
\label{tbl:Imag}
\begin{center}

{\normalsize

\begin{tabular}{|c||r|l||cc|c|}
\hline
 \(\text{cc}(M)\) & \(\lvert d\rvert\) & Factorization                           &   CT & \(\varkappa(K)\) &       \(\tau(K)\)  \\
\hline
   \(\mathbf{2}\) &         \(3\,896\) & \(=2^3\cdot 487\)                      & H.4  &       \((2111)\) & \((111,\mathbf{21},111,111)\) \\
   \(\mathbf{4}\) &        \(27\,156\) & \(=2^2\cdot 3\cdot 31\cdot 73\)        & F.11 &       \((4331)\) & \((111,\mathbf{32},32,111)\)  \\
   \(\mathbf{6}\) &       \(423\,640\) & \(=2^3\cdot 5\cdot 7\cdot 17\cdot 89\) & F.12 &       \((1422)\) & \((43,111,\mathbf{43},111)\)  \\
\hline
   \(\mathbf{8}\) &   \(99\,888\,340\) & \(=2^2\cdot 5\cdot 4994417\)           & F.13 &       \((3413)\) & \((111,\mathbf{54},111,65)\)  \\
\hline
\end{tabular}

}

\end{center}
\end{table}

%\newpage
%--------------------------------------------------------------------------------

\subsection{Real quadratic fields}
\label{ss:Real}

\noindent
The coclass \(\mathrm{cc}(M)\) of
\(M\simeq\mathrm{Gal}(\mathrm{F}_3^2(K)/K)\),
the second \(3\)-class group
of a real quadratic number field \(K\) with positive discriminant \(d>0\),
may be \textit{any} positive integer
\cite{Ma2012a,Ma2013,Ma2014a}.
Table \ref{tbl:Real} shows the minimal discriminants \(d\)
of real quadratic number fields \(K\) with assigned coclass \(\mathrm{cc}(M)\),
their prime factorization, the capitulation type (CT) \(\varkappa(K)\),
and the decisive logarithmic abelian type invariants, collected in the IPAD \(\tau(K)\).
Again the crucial \textbf{co-polarization} \(E_2\) with
\(\mathrm{lo}(\mathrm{Cl}_3(E_2))=\mathrm{cc}(M)+1\)
is pointed out in \textbf{bold} font.
For instance, in the fourth row, using the identification \(\mathbf{(32)}\,\hat{=}\,\lbrack 3^3,3^2\rbrack=\lbrack 27,9\rbrack\),
we have
\[
\text{cc}(M)+1=\mathbf{4}+1=5=3+2=\mathrm{lo}(\mathbf{(32)})=\mathrm{lo}(\lbrack 27,9\rbrack)=\mathrm{lo}(\mathrm{Cl}_3(E_2)).
\]

\renewcommand{\arraystretch}{1.1}

\begin{table}[ht]
\caption{Real quadratic number fields \(K\), OEIS A379524 \cite{OEIS2025}}
\label{tbl:Real}
\begin{center}

{\normalsize

\begin{tabular}{|c||r|l||cc|c|}
\hline
 \(\text{cc}(M)\) &             \(d\) & Factorization                &   CT & \(\varkappa(K)\) &       \(\tau(K)\)  \\
\hline
   \(\mathbf{1}\) &       \(32\,009\) & prime                        & a.3  &       \((0100)\) & \((\mathbf{11},21,11,11)\) \\
   \(\mathbf{2}\) &      \(214\,712\) & \(=2^3\cdot 26839\)          & G.19 &       \((2143)\) & \((\mathbf{21},21,21,21)\)  \\
   \(\mathbf{3}\) &      \(710\,652\) & \(=2^2\cdot 3\cdot 59221\)   & b.10 &       \((0320)\) & \((\mathbf{22},111,111,22)\)  \\
   \(\mathbf{4}\) &   \(8\,127\,208\) & \(=2^3\cdot 31\cdot 32771\)  & F.13 &       \((2443)\) & \((43,\mathbf{32},111,111)\)  \\
\hline
   \(\mathbf{5}\) & \(180\,527\,768\) & \(=2^3\cdot 569\cdot 39659\) & b.10 &       \((4001)\) & \((111,\mathbf{33},33,111)\)  \\
\hline
\end{tabular}

}

\end{center}
\end{table}

\noindent
The \(p\)-group generation algorithm \cite{HEO2005} by
Newman \cite{Nm1977} and
O'Brien \cite{Ob1990}
must be employed
in order to identify the groups \(M\) in Table
\ref{tbl:Imag} and \ref{tbl:Real}
exactly from their Artin pattern
\(\mathrm{AP}(M)=\left\lbrack\tau(K),\varkappa(K)\right\rbrack\),
where the capitulation type \(\varkappa(K)\) is definitely required,
too.

%\newpage
%--------------------------------------------------------------------------------

\subsection{Cyclic cubic fields}
\label{ss:Cycl}

\noindent
An extension to base fields of degree three is possible with the aid of
\cite{AoMa2024},
where cyclic cubic fields are constructed as
\(3\)-ray class fields modulo \(3\)-admissible conductors \(f\)
over the rational number field \(\mathbb{Q}\).
We computed all IPADs of cyclic cubic number fields \(K\)
with conductor \(f<10^6\) and
elementary bicyclic \(3\)-class group \(\mathrm{Cl}_3(K)\simeq (3,3)\).
This structure implies that the number of primes dividing the conductor \(f\)
is \(t\in\lbrace 2,3\rbrace\),
but the converse implication does not hold.
In significant contrast to quadratic base fields,
several non-isomorphic cyclic cubic fields,
collected in a multiplet \((K_\mu)_{\mu=1}^m\)
with multiplicity \(m=2^{t-1}\),
may share a common conductor \(f\).
Furthermore, not all members of a multiplet
share isomorphic \(3\)-class groups,
and we denote by \(\nu\le m\) the number of those
with elementary bicyclic \(3\)-class group
\(\mathrm{Cl}_3(K_\mu)\simeq (3,3)\).
Table \ref{tbl:Cycl} shows the minimal conductors \(f\)
of cyclic cubic number fields \(K\) with assigned coclass \(\mathrm{cc}(M)\)
of the second \(3\)-class group
\(M\simeq\mathrm{Gal}(\mathrm{F}_3^2(K)/K)\),
their prime factorization,
the number \(n\) and the multiplicity \(m\),
the corresponding graph in
\cite{AoMa2024},
the capitulation type (CT) \(\varkappa(K)\),
and the crucial IPAD \(\tau(K)\), containing the logarithmic abelian type invariants.
The decisive \textbf{co-polarization} \(E_2\) with
\(\mathrm{lo}(\mathrm{Cl}_3(E_2))=\mathrm{cc}(M)+1\)
is emphasized by \textbf{boldface} font.
As opposed to quadratic base fields,
\(M\) may be abelian with nilpotency class \(\text{cl}(M)=1\).

\renewcommand{\arraystretch}{1.1}

\begin{table}[ht]
\caption{Cyclic cubic number fields \(K\), OEIS A380103 \cite{OEIS2025}}
\label{tbl:Cycl}
\begin{center}

{\normalsize

\begin{tabular}{|cc||r|l|c|c||cc|c|}
\hline
 \(\text{cc}(M)\) & \(\text{cl}(M)\) &         \(f\) & Factorization           & \(\nu/m\) & Graph &   CT & \(\varkappa(K)\) &       \(\tau(K)\)  \\
\hline
   \(\mathbf{1}\) &            \(1\) &       \(657\) & \(=3^2\cdot 73\)        & \(2/2\) & 3     & a.1  &       \((0000)\) & \((\mathbf{1},1,1,1)\) \\
   \(\mathbf{1}\) &            \(2\) &    \(2\,439\) & \(=3^2\cdot 271\)       & \(2/2\) & 3     & A.1  &       \((1111)\) & \((\mathbf{11},2,2,2)\)  \\
   \(\mathbf{2}\) &                  &    \(7\,657\) & \(=13\cdot 19\cdot 31\) & \(3/4\) & I.2   & c.21 &       \((0231)\) & \((\mathbf{21},21,21,21)\)  \\
   \(\mathbf{3}\) &                  &   \(41\,839\) & \(=7\cdot 43\cdot 139\) & \(2/4\) & II.1  & b.10 &       \((0402)\) & \((\mathbf{22},111,22,111)\)  \\
   \(\mathbf{4}\) &                  &  \(231\,469\) & \(=7\cdot 43\cdot 769\) & \(4/4\) & III.8 & H.4  &       \((4442)\) & \((\mathbf{32},111,32,111)\)  \\
\hline
\end{tabular}

}

\end{center}
\end{table}

%\newpage
%--------------------------------------------------------------------------------

\subsection{Normal closures of pure cubic fields}
\label{ss:Meta}

\noindent
A further extension to base fields of degree six is based on
\cite{AoMa2025},
where the normal closures \(K=\mathbb{Q}(\zeta,\sqrt[3]{d})\) of pure cubic fields \(L=\mathbb{Q}(\sqrt[3]{d})\)
are constructed as \(3\)-ray class fields modulo \(3\)-admissible conductors \(f\)
over the quadratic cyclotomic field \(\mathbb{Q}(\zeta)\) with a primitive third root of unity \(\zeta\).
We computed all IPADs of pure metacyclic sextic number fields \(K\)
with conductor \(f<10^4\) and
elementary bicyclic \(3\)-class group \(\mathrm{Cl}_3(K)\simeq (3,3)\).
As opposed to quadratic base fields,
several non-isomorphic pure metacyclic sextic fields,
collected in a multiplet \((K_\mu)_{\mu=1}^m\)
with multiplicity \(m\), given by the formulas in
\cite[Thm. 2.1, p. 833]{Ma1992},
may share a common conductor \(f\).
Furthermore, not all members of a multiplet
share isomorphic second \(3\)-class groups,
and we denote by \(\nu\le m\) the number of those
with assigned fixed coclass \(\text{cc}(M)\)
and class \(\text{cl}(M)\).
Table \ref{tbl:Meta} shows the minimal conductors \(f\)
of pure metacyclic sextic number fields \(K\) with assigned coclass \(\mathrm{cc}(M)\)
of the second \(3\)-class group
\(M\simeq\mathrm{Gal}(\mathrm{F}_3^2(K)/K)\),
their prime factorization,
the number \(n\) and the multiplicity \(m\),
the Dedekind species,
the capitulation type (CT) \(\varkappa(K)\),
and the crucial logarithmic abelian type invariants contained in the IPAD \(\tau(K)\).
Again, the decisive \textbf{co-polarization} \(E_2\) with
\(\mathrm{lo}(\mathrm{Cl}_3(E_2))=\mathrm{cc}(M)+1\)
is emphasized by \textbf{boldface} font.
In contrast to quadratic base fields,
\(M\) may be abelian with nilpotency class \(\text{cl}(M)=1\).

\renewcommand{\arraystretch}{1.1}

\begin{table}[ht]
\caption{Pure metacyclic sextic number fields \(K\), OEIS A380104 \cite{OEIS2025}}
\label{tbl:Meta}
\begin{center}

{\normalsize

\begin{tabular}{|cc||r|l|c|c||cc|c|}
\hline
 \(\text{cc}(M)\) & \(\text{cl}(M)\) &         \(f\) & Factorization          & \(\nu/m\) & Species &   CT & \(\varkappa(K)\) &       \(\tau(K)\)  \\
\hline
   \(\mathbf{1}\) &            \(1\) &        \(30\) & \(=3\cdot 2\cdot 5\)   & \(1/1\) & IB      & a.1  &       \((0000)\) & \((\mathbf{1},1,1,1)\) \\
   \(\mathbf{1}\) &            \(2\) &        \(90\) & \(=3^2\cdot 2\cdot 5\) & \(2/4\) & IA      & A.1  &       \((1111)\) & \((\mathbf{11},2,2,2)\)  \\
   \(\mathbf{2}\) &                  &       \(418\) & \(=2\cdot 11\cdot 19\) & \(2/2\) & II      & b.10 &       \((4001)\) & \((111,\mathbf{21},21,111,)\)  \\
   \(\mathbf{3}\) &                  &    \(1\,626\) & \(=3\cdot 2\cdot 271\) & \(2/2\) & IB      & b.10 &       \((0402)\) & \((32,111,\mathbf{22},111)\)  \\
\hline
\end{tabular}

}

\end{center}
\end{table}

\noindent
Discriminants in terms of conductors are
\(d=f^2\) in Table
\ref{tbl:Cycl},
and \(d=-27\cdot f^4\) in Table
\ref{tbl:Meta}.

%\newpage
%--------------------------------------------------------------------------------

\section{Normal lattices of metabelian 3-groups}
\label{s:Normal}

\noindent
\textit{Class and coclass} can be \textit{visualized by the body of normal lattices}.
Since \(2\)-groups with elementary bicyclic abelianization \((2,2)\)
are CF-groups of maximal nilpotency class,
let \(p\) be an odd prime.
For \(p=3\), the following statements are proven by
\cite{Ne1989a},
for \(p\ge 5\) they may be but need not be true.
The lattice of all normal subgroups and, in particular,
the terms of the lower and upper central series
are determined for all metabelian \(p\)-groups \(G\) with generator rank \(d_1(G)=2\)
having abelianization of type \(G/G^\prime=(p,p)\) and minimal defect of commutativity \(k(G)=0\).
It is shown that many of these groups are realized as Galois groups of second Hilbert \(p\)-class fields
of an extensive set of quadratic base fields \(K=\mathbb{Q}(\sqrt{d})\),
which are characterized by capitulation types (CT) \(\varkappa(K)\) of \(p\)-classes.
Examples for cyclic cubic and pure metacyclic fields are also known.

%--------------------------------------------------------------------------------

Let \(p\ge 3\) be an odd prime number,
and \(G=\langle x,y\rangle\) be a two-generated metabelian \(p\)-group,
that is with soluble length \(\mathrm{sl}(G)\le 2\)
having an elementary bicyclic derived quotient \(G/G^\prime\) of type \((p,p)\).
Assume further that \(G\) is of order \(\#(G)=p^n\) with \(n\ge 2\),
and of nilpotency class \(c=\mathrm{cl}(G)\ge 1\).
Then \(G\) is of coclass \(r=\mathrm{cc}(G)=n-c\ge 1\).
Denote by
\[G=\gamma_1(G)>\gamma_2(G)=G^\prime>\ldots>\gamma_{c}(G)>\gamma_{c+1}(G)=1\]
the (descending) lower central series of \(G\), where \(\gamma_j(G)=\lbrack\gamma_{j-1}(G),G\rbrack\) for \(j\ge 2\), and by
\[1=\zeta_0(G)<\zeta_1(G)<\ldots<G^\prime=\zeta_{c-1}(G)<\zeta_{c}(G)=G\]
the (ascending) upper central series of \(G\), where \(\zeta_j(G)/\zeta_{j-1}(G)=\mathrm{Centre}(G/\zeta_{j-1}(G))\) for \(j\ge 1\).

Let \(s_2=t_2=\lbrack y,x\rbrack\) denote the main commutator of \(G\),
such that \(\gamma_2(G)=\langle s_2,\gamma_3(G)\rangle\) is cyclic of order \(p\).
By means of the two series
\(s_j=\lbrack s_{j-1},x\rbrack\) for \(j\ge 3\) and
\(t_\ell=\lbrack t_{\ell-1},y\rbrack\) for \(\ell\ge 3\)
of higher commutators and the subgroups
\(\Sigma_j=\langle s_j,\ldots,s_{c}\rangle\) with \(j\ge 3\) and \(T_\ell=\langle t_\ell,\ldots,t_{r+2}\rangle\) with \(\ell\ge 3\),
we obtain the following fundamental distinction of cases.

%--------------------------------------------------------------------------------

\begin{enumerate}
\item
The \textit{uniserial} case of a CF-group (\textit{cyclic factors}) of coclass \(\mathrm{cc}(G)=1\) (maximal class),
where \(t_3\in\Sigma_3\), \(\gamma_3(G)=\langle s_3,\gamma_4(G)\rangle\), \(r=1\), and \(c=n-1\).
There are two subcases:\\
(1.1) \(t_3=1\in\gamma_{c+1}(G)\), \(G\) contains an abelian maximal subgroup and \(k=k(G)=0\),\\
(1.2) \(1\ne t_3\in\gamma_{c+1-k}(G)\), with positive defect of commutativity \(1\le k=k(G)\le\min(c-3,p-2)\),
where all maximal subgroups are non-abelian.
\item
The \textit{biserial} case of a non-CF- or BCF-group (\textit{bicyclic or cyclic factors}) of coclass \(\mathrm{cc}(G)\ge 2\),
where \(t_3\not\in\Sigma_3\), \(\gamma_3(G)=\langle s_3,t_3,\gamma_4(G)\rangle\), \(r\ge 2\), and \(c\le n-2\).
Again there exist two subcases, characterized by the \textit{defect of commutativity} \(k(G)\) of \(G\):\\
(2.1) \(t_{r+2}=1\in\gamma_{c+1}(G)\), where \(\Sigma_3\cap T_3=1\) and \(k=k(G)=0\),\\
(2.2) \(1\ne t_{r+2}\in\gamma_{c+1-k}(G)\), for some \(k=k(G)\ge 1\), where \(\Sigma_3\cap T_3\le\gamma_{c+1-k}(G)\).
\end{enumerate}

%--------------------------------------------------------------------------------

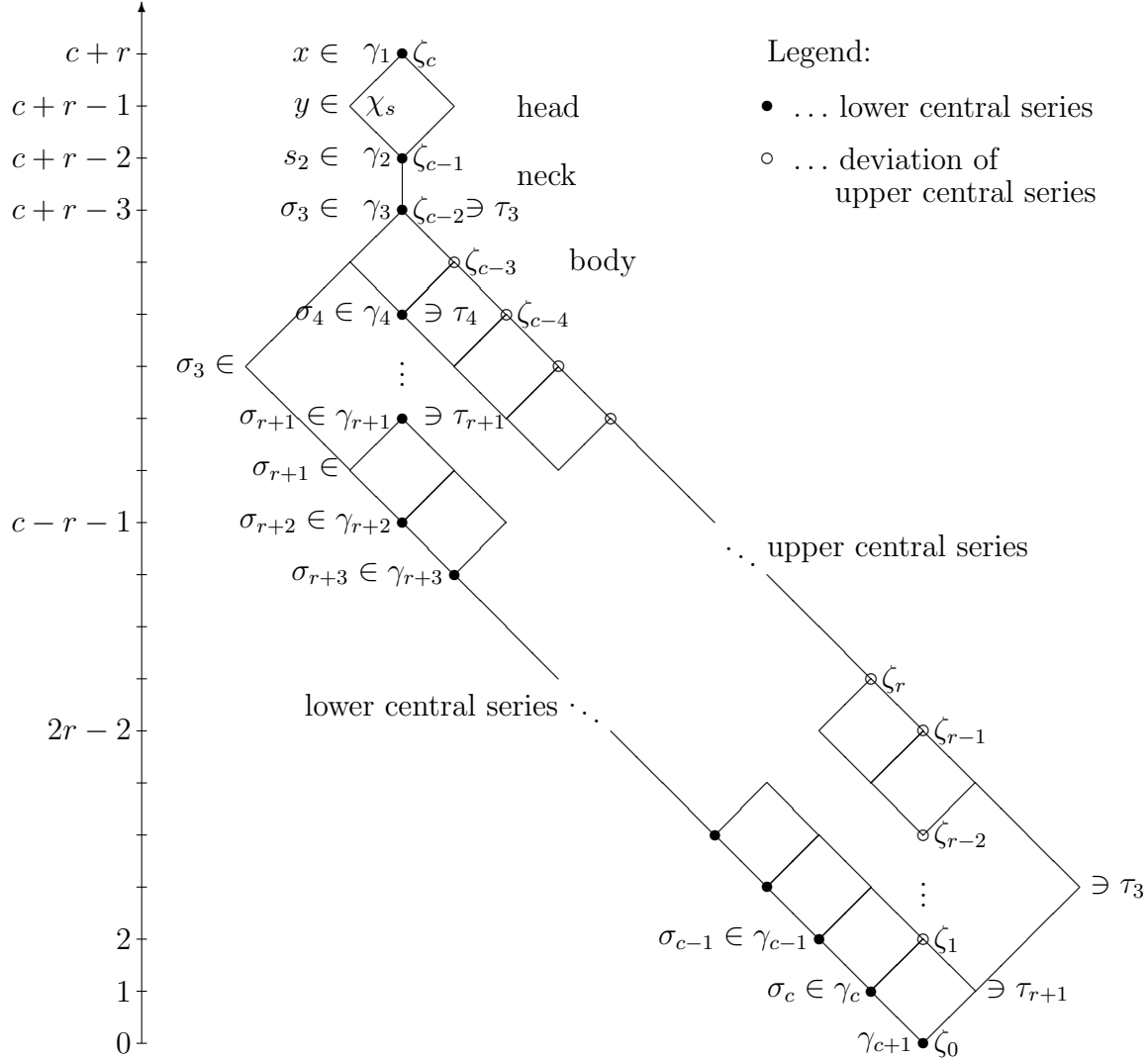
\begin{figure}[ht]
\caption{Full normal lattice, including lower and upper central series, of a \(p\)-group \(G\) with \(G/G^\prime\simeq (p,p)\), \(\mathrm{cl}(G)=c\), \(\mathrm{cc}(G)=r\), \(\mathrm{sl}(G)=2\), \(k(G)=0\).}
\label{fig:NormalLatticeFigure}

% Full normal lattice, including lower and upper central series

\setlength{\unitlength}{0.7cm}
\begin{picture}(18,22)(-5,-1)

% scale of orders
\put(-5,20.3){\makebox(0,0)[cb]{Logarithmic Order}}
\put(-5,19){\vector(0,1){1}}
\put(-5.2,19){\makebox(0,0)[rc]{\({c+r}\)}}
\put(-5.2,18){\makebox(0,0)[rc]{\({c+r-1}\)}}
\put(-5.2,17){\makebox(0,0)[rc]{\({c+r-2}\)}}
\put(-5.2,16){\makebox(0,0)[rc]{\({c+r-3}\)}}
\put(-5.2,10){\makebox(0,0)[rc]{\({c-r-1}\)}}
\put(-5.2,6){\makebox(0,0)[rc]{\({2r-2}\)}}
\put(-5.2,2){\makebox(0,0)[rc]{\(2\)}}
\put(-5.2,1){\makebox(0,0)[rc]{\(1\)}}
\put(-5.2,0){\makebox(0,0)[rc]{\(0\)}}
\multiput(-5.1,0)(0,1){20}{\line(1,0){0.2}}
\put(-5,0){\line(0,1){19}}

% Legend
\put(7,19){\makebox(0,0)[lc]{Legend:}}
\put(7,18){\circle*{0.2}}
\put(7.5,18){\makebox(0,0)[lc]{\(\ldots\) lower central series}}
\put(7,17){\circle{0.2}}
\put(7.5,17){\makebox(0,0)[lc]{\(\ldots\) deviation of}}
\put(8.3,16.3){\makebox(0,0)[lc]{upper central series}}

% full normal lattice

% top of head
\put(-0.2,19){\makebox(0,0)[rc]{\(\gamma_1\)}}
\put(0,19){\circle*{0.2}}
\put(0.2,19){\makebox(0,0)[lc]{\(\zeta_{c}\)}}

% head
\put(-1.2,19){\makebox(0,0)[rc]{\(x\in\)}}
\multiput(0,19)(-1,-1){2}{\line(1,-1){1}}
\multiput(0,19)(1,-1){2}{\line(-1,-1){1}}
\put(-0.7,18){\makebox(0,0)[lc]{\(\chi_s\)}}
\put(-1.2,18){\makebox(0,0)[rc]{\(y\in\)}}

\put(2.2,18){\makebox(0,0)[lc]{head}}

% bottom of head
\put(-1.2,17){\makebox(0,0)[rc]{\(s_2\in\)}}
\put(-0.2,17){\makebox(0,0)[rc]{\(\gamma_2\)}}
\put(0,17){\circle*{0.2}}
\put(0.2,17){\makebox(0,0)[lc]{\(\zeta_{c-1}\)}}

% (bottle) neck
\put(0,17){\line(0,-1){1}}

\put(2.2,16.7){\makebox(0,0)[lc]{neck}}

% top of body
\put(-1.2,16){\makebox(0,0)[rc]{\(\sigma_3\in\)}}
\put(-0.2,16){\makebox(0,0)[rc]{\(\gamma_3\)}}
\put(1.2,16){\makebox(0,0)[lc]{\(\ni\tau_3\)}}
\put(0,16){\circle*{0.2}}
\put(0.2,16){\makebox(0,0)[lc]{\(\zeta_{c-2}\)}}

% body
\multiput(0,16)(-1,-1){2}{\line(1,-1){1}}
\multiput(0,16)(1,-1){2}{\line(-1,-1){1}}

\put(-0.2,14){\makebox(0,0)[rc]{\(\sigma_4\in\gamma_4\)}}
\put(0.4,14){\makebox(0,0)[lc]{\(\ni\tau_4\)}}
\put(0,14){\circle*{0.2}}

\put(1,15){\circle{0.2}}
\put(1.2,15){\makebox(0,0)[lc]{\(\zeta_{c-3}\)}}

\put(3.2,15){\makebox(0,0)[lc]{body}}

\multiput(1,15)(-1,-1){2}{\line(1,-1){1}}
\multiput(1,15)(1,-1){2}{\line(-1,-1){1}}

\put(2,14){\circle{0.2}}
\put(2.2,14){\makebox(0,0)[lc]{\(\zeta_{c-4}\)}}

\multiput(2,14)(-1,-1){2}{\line(1,-1){1}}
\multiput(2,14)(1,-1){2}{\line(-1,-1){1}}

\put(3,13){\circle{0.2}}

\multiput(3,13)(-1,-1){2}{\line(1,-1){1}}
\multiput(3,13)(1,-1){2}{\line(-1,-1){1}}

\put(4,12){\circle{0.2}}

\put(-3.2,13){\makebox(0,0)[rc]{\(\sigma_3\in\)}}

\put(4,12){\line(1,-1){2}}

% upper left corner of rectangle
\put(-3,13){\line(1,1){2}}
\put(-3,13){\line(1,-1){2}}

\put(0,13){\makebox(0,0)[cc]{\(\vdots\)}}

\put(-1.2,11){\makebox(0,0)[rc]{\(\sigma_{r+1}\in\)}}

\put(-0.2,12){\makebox(0,0)[rc]{\(\sigma_{r+1}\in\gamma_{r+1}\)}}
\put(0.4,12){\makebox(0,0)[lc]{\(\ni\tau_{r+1}\)}}
\put(0,12){\circle*{0.2}}

\multiput(0,12)(-1,-1){2}{\line(1,-1){1}}
\multiput(0,12)(1,-1){2}{\line(-1,-1){1}}

\put(-0.2,10){\makebox(0,0)[rc]{\(\sigma_{r+2}\in\gamma_{r+2}\)}}
\put(0,10){\circle*{0.2}}

\multiput(1,11)(-1,-1){2}{\line(1,-1){1}}
\multiput(1,11)(1,-1){2}{\line(-1,-1){1}}

\put(0.8,9){\makebox(0,0)[rc]{\(\sigma_{r+3}\in\gamma_{r+3}\)}}
\put(1,9){\circle*{0.2}}

\put(1,9){\line(1,-1){2}}
\put(3.5,6.5){\makebox(0,0)[cc]{\(\ddots\)}}
\put(3,6.5){\makebox(0,0)[rc]{lower central series}}

\put(7,9.5){\makebox(0,0)[lc]{upper central series}}
\put(6.5,9.5){\makebox(0,0)[cc]{\(\ddots\)}}
\put(9,7){\line(-1,1){2}}

\put(9,7){\circle{0.2}}
\put(9.2,7){\makebox(0,0)[lc]{\(\zeta_{r}\)}}

\multiput(9,7)(-1,-1){2}{\line(1,-1){1}}
\multiput(9,7)(1,-1){2}{\line(-1,-1){1}}

\put(10,6){\circle{0.2}}
\put(10.2,6){\makebox(0,0)[lc]{\(\zeta_{r-1}\)}}

\multiput(10,6)(-1,-1){2}{\line(1,-1){1}}
\multiput(10,6)(1,-1){2}{\line(-1,-1){1}}

\put(10,4){\circle{0.2}}
\put(10.2,4){\makebox(0,0)[lc]{\(\zeta_{r-2}\)}}
\put(13.2,3){\makebox(0,0)[lc]{\(\ni\tau_3\)}}

% lower right corner of rectangle
\put(13,3){\line(-1,1){2}}
\put(13,3){\line(-1,-1){2}}

\put(10,3){\makebox(0,0)[cc]{\(\vdots\)}}

\put(6,4){\line(-1,1){2}}

\put(6,4){\circle*{0.2}}

\multiput(7,5)(-1,-1){2}{\line(1,-1){1}}
\multiput(7,5)(1,-1){2}{\line(-1,-1){1}}

\put(7,3){\circle*{0.2}}

\multiput(8,4)(-1,-1){2}{\line(1,-1){1}}
\multiput(8,4)(1,-1){2}{\line(-1,-1){1}}

\put(7.8,2){\makebox(0,0)[rc]{\(\sigma_{c-1}\in\gamma_{c-1}\)}}
\put(8,2){\circle*{0.2}}

\multiput(9,3)(-1,-1){2}{\line(1,-1){1}}
\multiput(9,3)(1,-1){2}{\line(-1,-1){1}}

\put(8.8,1){\makebox(0,0)[rc]{\(\sigma_{c}\in\gamma_{c}\)}}
\put(9,1){\circle*{0.2}}

\put(10,2){\circle{0.2}}
\put(10.2,2){\makebox(0,0)[lc]{\(\zeta_1\)}}
\put(11.2,1){\makebox(0,0)[lc]{\(\ni\tau_{r+1}\)}}

% tail
\multiput(10,2)(-1,-1){2}{\line(1,-1){1}}
\multiput(10,2)(1,-1){2}{\line(-1,-1){1}}

% bottom of tail
\put(9.8,0){\makebox(0,0)[rc]{\(\gamma_{c+1}\)}}
\put(10,0){\circle*{0.2}}
\put(10.2,0){\makebox(0,0)[lc]{\(\zeta_0\)}}

\end{picture}

\end{figure}

In this article, we are interested in two-generator metabelian \(p\)-groups \(G=\langle x,y\rangle\)
of coclass \(\mathrm{cc}(G)\ge 2\) having the convenient property \(\Sigma_3\cap T_3=1\), respectively \(k(G)=0\),
where the product \(\Sigma_3\times T_3\) is direct and coincides with the major part of the \textit{normal lattice} of \(G\),
as shown in Figure
\ref{fig:NormalLatticeFigure},
with \(s_j,t_j\) replaced by \(\sigma_j,\tau_j\), for \(j\ge 3\), as in section \S\
\ref{s:HeteroHomo}.

\begin{definition}
\label{dfn:Diamond}
A pair \((U,V)\) of normal subgroups of a \(p\)-group \(G\),
such that \(V<U\le G\) and \((U:V)=p^2\),
is called a \textit{diamond} if the quotient \(U/V\) is abelian of type \((p,p)\).
\end{definition}

If \((U,V)\) is a diamond and \(U=\langle u_1,u_2,V\rangle\),
then the \(p+1\) intermediate subgroups of \(G\) between \(U\) and \(V\) are given by
\(\langle u_2,V\rangle\) and \(\langle u_1u_2^{i-2},V\rangle\) with \(2\le i\le p+1\).

%--------------------------------------------------------------------------------

\medskip
%\noindent
In this section, let \(G=\langle x,y\rangle\) be a metabelian \(p\)-group
with two generators \(x,y\), having abelianization \(G/G^\prime\) of type \((p,p)\)
and satisfying the independence condition \(\Sigma_3\cap T_3=1\), that is,
\(G\) is a metabelian \(p\)-group with defect of commutativity \(k(G)=0\)
\cite[\S\ 3.1.1, p. 412, and \S\ 3.3.2, p. 429]{Ma2013}.
We assume that \(G\) is of coclass \(\mathrm{cc}(G)\ge 2\),
since the normal lattice of \(p\)-groups of maximal class
has been determined by Blackburn
\cite{Bl1958}.

%--------------------------------------------------------------------------------

\begin{theorem}
\label{thm:NormalLattice}
The complete normal lattice of \(G\)
%is given by
contains
the heading diamond \((G,G^\prime)\)
and the rectangle \(\bigl((P_{j,\ell},P_{j+1,\ell+1})\bigr)_{3\le j\le c,\ 3\le\ell\le r+1}\) of trailing diamonds,
where \(P_{j,\ell}=\Sigma_j\times T_\ell\) for \(3\le j\le c+1\), \(3\le\ell\le r+2\).
The structure of the normal lattice is visualized in Figure
\ref{fig:NormalLatticeFigure}.
\end{theorem}

%\noindent
Note that
\(P_{j,\ell}=\langle s_j,\ldots,s_{c}\rangle\times\langle t_\ell,\ldots,t_{r+1}\rangle=\langle s_j,t_\ell,P_{j+1,\ell+1}\rangle\)
for \(3\le j\le c\), \(3\le\ell\le r+1\).

\begin{corollary}
\label{cor:NormalLatticeComplete}
The complete normal lattice of \(G\)
consists exactly of the normal subgroups given in Theorem
\ref{thm:NormalLattice}.
\end{corollary}

%--------------------------------------------------------------------------------

\begin{corollary}
\label{cor:NormalLatticeCount}
The total number of normal subgroups of \(G\) is given by
\[1+c+r+p\cdot\lbrack 3+cr-c-2r\rbrack,\]
in particular, for \(p=3\) it is given by
\[10+3cr-2c-5r.\]
\end{corollary}

%--------------------------------------------------------------------------------

\begin{corollary}
\label{cor:CentralSeries}
Blackburn's two-step centralizers of \(G\)
\cite{Bl1958}
are given by

\begin{equation*}
\chi_j(G)=
\begin{cases}
G^\prime                  \text{ for } 1\le j\le r,\\
\langle y,G^\prime\rangle \text{ for } r+1\le j\le c-1,\\
G                         \text{ for } j\ge c.
\end{cases}
\end{equation*}

\noindent
None of the maximal subgroups of \(G\)
occurs as a two-step centralizer, when \(c=r+1\).

\begin{enumerate}

\item
The factors of the lower central series of \(G\) are given by

\begin{equation*}
\gamma_j(G)/\gamma_{j+1}(G)\simeq
\begin{cases}
(p,p) \text{ for } j=1 \text{ and } 3\le j\le r+1,\\
(p)   \text{ for } j=2 \text{ and } r+2\le j\le c.
\end{cases}
\end{equation*}

\item
The terms of the lower central series of \(G\) are given by

\begin{equation*}
\gamma_j(G)=
\begin{cases}
\langle x,y,G^\prime\rangle    \text{ for } j=1,\\
\langle s_2,\gamma_3(G)\rangle \text{ for } j=2,\\
P_{j,j}                        \text{ for } 3\le j\le r+1,\\
\Sigma_j                       \text{ for } r+2\le j\le c.
\end{cases}
\end{equation*}

\item
The factors of the upper central series of \(G\) are given by

\begin{equation*}
\zeta_j(G)/\zeta_{j-1}(G)\simeq
\begin{cases}
(p,p) \text{ for } 1\le j\le r-1 \text{ and } j=c,\\
(p)   \text{ for } r\le j\le c-1.
\end{cases}
\end{equation*}

\item
The terms of the upper central series of \(G\) are given by

\begin{equation*}
\zeta_j(G)=
\begin{cases}
P_{c+1-j,r+2-j}                     \text{ for } 1\le j\le r-1,\\
P_{c+1-j,3}                         \text{ for } r\le j\le c-2,\\
\langle s_2,\zeta_{c-2}(G)\rangle \text{ for } j=c-1,\\
\langle x,y,\zeta_{c-1}(G)\rangle \text{ for } j=c.
\end{cases}
\end{equation*}

\end{enumerate}

\end{corollary}

\begin{proof}
The heading diamond \((G,G^\prime)\) has 
a top vertex, a bottom vertex, and \(p+1\) intermediate vertices,
that is, a total of \(p+3\) vertices.
The rectangle \(\bigl((P_{j,\ell},P_{j+1,\ell+1})\bigr)_{3\le j\le c,\ 3\le\ell\le r+1}\)
consists of \((c-2)\cdot (r-1)\) trailing diamonds, each with \(p-1\) inner vertices,
and a total of \((c-1)\cdot r\) outer vertices.
The number of all normal subgroups is therefore given by
\begin{equation*}
\begin{aligned}
    & p+3+(c-2)(r-1)(p-1)+(c-1)r \\
=\, & p+3+(cr-c-2r+2)(p-1)+cr-r \\
=\, & 3+p(cr-c-2r+3)-cr+c+2r-2+cr-r \\
=\, & 1+c+r+p(3+cr-c-2r)
\end{aligned}
\end{equation*}
generally, and for \(p=3\) by \(1+c+r+9+3cr-3c-6r=10+3cr-2c-5r\).
This proves Corollary
\ref{cor:NormalLatticeCount}.

The remaining claims of Theorem
\ref{thm:NormalLattice}
and the Corollaries
\ref{cor:NormalLatticeComplete}
and
\ref{cor:CentralSeries}
for a group \(G\) with defect \(k(G)=0\)
(and thus a bicyclic center \(\zeta_1(G)\)
and a direct product \(P_{3,3}=\Sigma_3\times T_3\))
are consequences of the definition of the subgroups
\(P_{j,\ell}=\Sigma_j\times T_\ell\) for \(3\le j\le c+1\), \(3\le\ell\le r+2\),
and the results in Nebelung's Ph.D. thesis \cite{Ne1989a},
which are summarized in
\cite[\S\ 8, pp. 456--461]{Ma2014a}.
\end{proof}

%\newpage
%--------------------------------------------------------------------------------

%\newpage
%--------------------------------------------------------------------------------

\subsection{Cyclic factor groups}
\label{ss:CF}

\noindent
CF-groups \(G\) have maximal nilpotency class (coclass \(r=1\))
with \textit{cyclic factors} of the central series,
except the head \(\gamma_1(G)/\gamma_2(G)\).
See Figure
\ref{fig:CFGroups}.
The rectangle of trailing diamonds degenerates to a line, and
the upper central series is the reverse lower central series.
As second \(p\)-class groups \(M=\mathrm{Gal}(\mathrm{F}_p^2(K)\vert K)\)
they cannot be realized by imaginary quadratic fields \(K\),
and by \textit{real quadratic} fields \(K=\mathbb{Q}(\sqrt{d})\)
only for \textit{odd} nilpotency class \(c\ge 3\),
since \(M\) can never be abelian (\(c=1\)) for a quadratic field \(K\).

\begin{example}
\label{exm:CFGroups}
CF-groups \(M\) with \(p\in\lbrace 3,5,7\rbrace\), specified by capitulation types (CT). \\
\(c=3\) for \(d\in\lbrace 32\,009, 72\,329, 142\,097\rbrace\) with CT a.3, a.2, a.3\({}^\ast\)
\cite[Fig. 2, p. 143]{Ma2017}, \\
\(c=5\) for \(d\in\lbrace 62\,501, 152\,949, 252\,977, 494\,236, 790\,085\rbrace\), 3 a.1, a.3, a.2
\cite[Fig. 1, p. 44]{Ma2018a}, \\
\(c=7\) for \(d\in\lbrace 2\,905\,160, 10\,200\,108, 10\,399\,596, 14\,458\,876, 27\,780\,297\rbrace\), excited states,\\
\(c=9\) for \(d\in\lbrace 37\,304\,664, 40\,980\,808, 62\,565\,429, 63\,407\,037, 208\,540\,653\rbrace\),
if \(p=3\), \\
but the same \textit{selection rule} also holds for bigger prime numbers \(p>3\): \\
\(c=3\) for \(d\in\lbrace 244\,641, 1\,167\,541\rbrace\), \\
\(c=5\) for \(d\in\lbrace 1\,129\,841, 3\,812\,377\rbrace\),
if \(p=5\)
\cite[Fig. 3.3, p. 425, Tbl. 3.4, p. 427]{Ma2013},
and \\
\(c=3\) for \(d\in\lbrace 1\,633\,285, 2\,713\,121\rbrace\),
if \(p=7\)
\cite[Fig. 3.4, Tbl. 3.5, p. 428]{Ma2013}. \\
For \textit{cyclic cubic} fields \(K\) with conductor \(f\) and \(p=3\),
there are no constraints for the parity of the nilpotency class by selection rules, e.g. (Table
\ref{tbl:Cycl} and
\cite{AoMa2024}) \\
\(c=1\) for \(f=657\), \(c=2\) for \(f=2439\), \(c=3\) for \(f\in\lbrace 3913, 8001\rbrace\),  and \(c=4\) for \(f=4977\). \\
Similarly for \textit{pure metacyclic sextic} fields \(K\) with conductor \(f\) and \(p=3\), e.g. \\
\(c=1\) for \(f=30\), \(c=2\) for \(f=90\) (Table
\ref{tbl:Meta}).
(Even class \(c\) is not drawn in Figure \ref{fig:CFGroups}.)
\end{example}

\begin{figure}[ht]
\caption{CF-groups \(M=\mathrm{Gal}(\mathrm{F}_p^2(K)\vert K)\) with cyclic factors only.}
\label{fig:CFGroups}

{\small

% CF Groups

\setlength{\unitlength}{0.5cm}
\begin{picture}(29,12)(-7,0)

% scale of orders
%\put(-5.8,14.3){\makebox(0,0)[cb]{order \(3^n\), \(p^n\)}}
%\put(-5,13){\vector(0,1){1}}
\put(-5.8,11.3){\makebox(0,0)[cb]{Order \(3^n\), \(p^n\)}}
\put(-5,10){\vector(0,1){1}}
%\put(-5.2,13){\makebox(0,0)[rc]{\(1594323\)}}
%\put(-4.8,13){\makebox(0,0)[lc]{\(p^{13}\)}}
%\put(-5.2,12){\makebox(0,0)[rc]{\(531441\)}}
%\put(-4.8,12){\makebox(0,0)[lc]{\(p^{12}\)}}
%\put(-5.2,11){\makebox(0,0)[rc]{\(177147\)}}
%\put(-4.8,11){\makebox(0,0)[lc]{\(p^{11}\)}}
\put(-5.2,10){\makebox(0,0)[rc]{\(59049\)}}
\put(-4.8,10){\makebox(0,0)[lc]{\(p^{10}\)}}
\put(-5.2,9){\makebox(0,0)[rc]{\(19683\)}}
\put(-4.8,9){\makebox(0,0)[lc]{\(p^9\)}}
\put(-5.2,8){\makebox(0,0)[rc]{\(6561\)}}
\put(-4.8,8){\makebox(0,0)[lc]{\(p^8\)}}
\put(-5.2,7){\makebox(0,0)[rc]{\(2187\)}}
\put(-4.8,7){\makebox(0,0)[lc]{\(p^7\)}}
\put(-5.2,6){\makebox(0,0)[rc]{\(729\)}}
\put(-4.8,6){\makebox(0,0)[lc]{\(p^6\)}}
\put(-5.2,5){\makebox(0,0)[rc]{\(243\)}}
\put(-4.8,5){\makebox(0,0)[lc]{\(p^5\)}}
\put(-5.2,4){\makebox(0,0)[rc]{\(81\)}}
\put(-4.8,4){\makebox(0,0)[lc]{\(p^4\)}}
\put(-5.2,3){\makebox(0,0)[rc]{\(27\)}}
\put(-4.8,3){\makebox(0,0)[lc]{\(p^3\)}}
\put(-5.2,2){\makebox(0,0)[rc]{\(9\)}}
\put(-4.8,2){\makebox(0,0)[lc]{\(p^2\)}}
\put(-5.2,1){\makebox(0,0)[rc]{\(3\)}}
\put(-4.8,1){\makebox(0,0)[lc]{\(p\)}}
\put(-5.2,0){\makebox(0,0)[rc]{\(1\)}}
\multiput(-5.1,0)(0,1){11}{\line(1,0){0.2}}
\put(-5,0){\line(0,1){10}}
%\multiput(-5.1,0)(0,1){14}{\line(1,0){0.2}}
%\put(-5,0){\line(0,1){13}}

% e = 2, m = 10
\put(1,10.5){\makebox(0,0)[cb]{\(r=1\), \(c=9\)}}
\put(1,10){\circle*{0.2}}
\multiput(1,10)(-1,-1){2}{\line(1,-1){1}}
\multiput(1,10)(1,-1){2}{\line(-1,-1){1}}
% bottle neck and body
\multiput(1,8)(0,-1){9}{\circle*{0.2}}
\multiput(1,8)(0,-1){8}{\line(0,-1){1}}

% e = 2, m = 8
\put(6,8.5){\makebox(0,0)[cb]{\(r=1\), \(c=7\)}}
\put(6,8){\circle*{0.2}}
\multiput(6,8)(-1,-1){2}{\line(1,-1){1}}
\multiput(6,8)(1,-1){2}{\line(-1,-1){1}}
% bottle neck and body
\multiput(6,6)(0,-1){7}{\circle*{0.2}}
\multiput(6,6)(0,-1){6}{\line(0,-1){1}}

% e = 2, m = 6
\put(11,6.5){\makebox(0,0)[cb]{\(r=1\), \(c=5\)}}
\put(11,6){\circle*{0.2}}
\multiput(11,6)(-1,-1){2}{\line(1,-1){1}}
\multiput(11,6)(1,-1){2}{\line(-1,-1){1}}
% bottle neck and body
\multiput(11,4)(0,-1){5}{\circle*{0.2}}
\multiput(11,4)(0,-1){4}{\line(0,-1){1}}

% e = 2, m = 4
\put(16,4.5){\makebox(0,0)[cb]{\(r=1\), \(c=3\)}}
\put(16,4){\circle*{0.2}}
\multiput(16,4)(-1,-1){2}{\line(1,-1){1}}
\multiput(16,4)(1,-1){2}{\line(-1,-1){1}}
% bottle neck and body
\multiput(16,2)(0,-1){3}{\circle*{0.2}}
\multiput(16,2)(0,-1){2}{\line(0,-1){1}}

% e = 2, m = 2
\put(21,2.5){\makebox(0,0)[cb]{\(r=1\), \(c=1\)}}
\put(21,2){\circle*{0.2}}
\multiput(21,2)(-1,-1){2}{\line(1,-1){1}}
\multiput(21,2)(1,-1){2}{\line(-1,-1){1}}
\multiput(21,0)(0,-1){1}{\circle*{0.2}}

\end{picture}

}

\end{figure}

%\newpage
%--------------------------------------------------------------------------------

\subsection{Bicyclic factor groups}
\label{ss:BF}

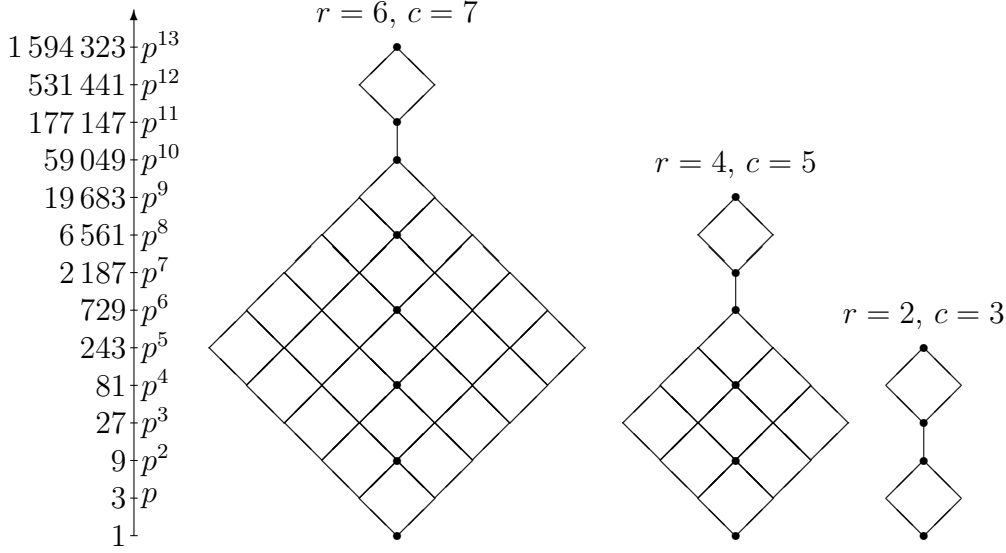
\begin{figure}[ht]
\caption{BF-groups \(M=\mathrm{Gal}(\mathrm{F}_p^2(K)\vert K)\) with bicyclic factors only.}
\label{fig:BFGroups}

% BF Groups

\setlength{\unitlength}{0.5cm}
\begin{picture}(24,15)(-7,0)

% scale of orders
\put(-5.8,14.3){\makebox(0,0)[cb]{Order \(3^n\), \(p^n\)}}
\put(-5,13){\vector(0,1){1}}
\put(-5.2,13){\makebox(0,0)[rc]{\(1\,594\,323\)}}
\put(-4.8,13){\makebox(0,0)[lc]{\(p^{13}\)}}
\put(-5.2,12){\makebox(0,0)[rc]{\(531\,441\)}}
\put(-4.8,12){\makebox(0,0)[lc]{\(p^{12}\)}}
\put(-5.2,11){\makebox(0,0)[rc]{\(177\,147\)}}
\put(-4.8,11){\makebox(0,0)[lc]{\(p^{11}\)}}
\put(-5.2,10){\makebox(0,0)[rc]{\(59\,049\)}}
\put(-4.8,10){\makebox(0,0)[lc]{\(p^{10}\)}}
\put(-5.2,9){\makebox(0,0)[rc]{\(19\,683\)}}
\put(-4.8,9){\makebox(0,0)[lc]{\(p^9\)}}
\put(-5.2,8){\makebox(0,0)[rc]{\(6\,561\)}}
\put(-4.8,8){\makebox(0,0)[lc]{\(p^8\)}}
\put(-5.2,7){\makebox(0,0)[rc]{\(2\,187\)}}
\put(-4.8,7){\makebox(0,0)[lc]{\(p^7\)}}
\put(-5.2,6){\makebox(0,0)[rc]{\(729\)}}
\put(-4.8,6){\makebox(0,0)[lc]{\(p^6\)}}
\put(-5.2,5){\makebox(0,0)[rc]{\(243\)}}
\put(-4.8,5){\makebox(0,0)[lc]{\(p^5\)}}
\put(-5.2,4){\makebox(0,0)[rc]{\(81\)}}
\put(-4.8,4){\makebox(0,0)[lc]{\(p^4\)}}
\put(-5.2,3){\makebox(0,0)[rc]{\(27\)}}
\put(-4.8,3){\makebox(0,0)[lc]{\(p^3\)}}
\put(-5.2,2){\makebox(0,0)[rc]{\(9\)}}
\put(-4.8,2){\makebox(0,0)[lc]{\(p^2\)}}
\put(-5.2,1){\makebox(0,0)[rc]{\(3\)}}
\put(-4.8,1){\makebox(0,0)[lc]{\(p\)}}
\put(-5.2,0){\makebox(0,0)[rc]{\(1\)}}
\multiput(-5.1,0)(0,1){14}{\line(1,0){0.2}}
\put(-5,0){\line(0,1){13}}

% e = 7, m = 8
\put(2,13.5){\makebox(0,0)[cb]{\(r=6\), \(c=7\)}}
\put(2,13){\circle*{0.2}}
\multiput(2,13)(-1,-1){2}{\line(1,-1){1}}
\multiput(2,13)(1,-1){2}{\line(-1,-1){1}}
\put(2,11){\circle*{0.2}}
% bottle neck
\put(2,11){\line(0,-1){1}}
% body
\put(2,10){\circle*{0.2}}
\multiput(2,10)(-1,-1){2}{\line(1,-1){1}}
\multiput(2,10)(1,-1){2}{\line(-1,-1){1}}
\put(2,8){\circle*{0.2}}
\multiput(2,8)(-1,-1){2}{\line(1,-1){1}}
\multiput(2,8)(1,-1){2}{\line(-1,-1){1}}
\put(2,6){\circle*{0.2}}
\multiput(2,6)(-1,-1){2}{\line(1,-1){1}}
\multiput(2,6)(1,-1){2}{\line(-1,-1){1}}
\put(2,4){\circle*{0.2}}
\multiput(2,4)(-1,-1){2}{\line(1,-1){1}}
\multiput(2,4)(1,-1){2}{\line(-1,-1){1}}
\put(2,2){\circle*{0.2}}
\multiput(2,2)(-1,-1){2}{\line(1,-1){1}}
\multiput(2,2)(1,-1){2}{\line(-1,-1){1}}
\put(2,0){\circle*{0.2}}

\multiput(1,9)(-1,-1){2}{\line(1,-1){1}}
\multiput(1,9)(1,-1){2}{\line(-1,-1){1}}
\multiput(3,9)(-1,-1){2}{\line(1,-1){1}}
\multiput(3,9)(1,-1){2}{\line(-1,-1){1}}

\multiput(0,8)(-1,-1){2}{\line(1,-1){1}}
\multiput(0,8)(1,-1){2}{\line(-1,-1){1}}
\multiput(4,8)(-1,-1){2}{\line(1,-1){1}}
\multiput(4,8)(1,-1){2}{\line(-1,-1){1}}

\multiput(-1,7)(-1,-1){2}{\line(1,-1){1}}
\multiput(-1,7)(1,-1){2}{\line(-1,-1){1}}
\multiput(1,7)(-1,-1){2}{\line(1,-1){1}}
\multiput(1,7)(1,-1){2}{\line(-1,-1){1}}
\multiput(3,7)(-1,-1){2}{\line(1,-1){1}}
\multiput(3,7)(1,-1){2}{\line(-1,-1){1}}
\multiput(5,7)(-1,-1){2}{\line(1,-1){1}}
\multiput(5,7)(1,-1){2}{\line(-1,-1){1}}

\multiput(-2,6)(-1,-1){2}{\line(1,-1){1}}
\multiput(-2,6)(1,-1){2}{\line(-1,-1){1}}
\multiput(0,6)(-1,-1){2}{\line(1,-1){1}}
\multiput(0,6)(1,-1){2}{\line(-1,-1){1}}
\multiput(4,6)(-1,-1){2}{\line(1,-1){1}}
\multiput(4,6)(1,-1){2}{\line(-1,-1){1}}
\multiput(6,6)(-1,-1){2}{\line(1,-1){1}}
\multiput(6,6)(1,-1){2}{\line(-1,-1){1}}

\multiput(-1,5)(-1,-1){2}{\line(1,-1){1}}
\multiput(-1,5)(1,-1){2}{\line(-1,-1){1}}
\multiput(1,5)(-1,-1){2}{\line(1,-1){1}}
\multiput(1,5)(1,-1){2}{\line(-1,-1){1}}
\multiput(3,5)(-1,-1){2}{\line(1,-1){1}}
\multiput(3,5)(1,-1){2}{\line(-1,-1){1}}
\multiput(5,5)(-1,-1){2}{\line(1,-1){1}}
\multiput(5,5)(1,-1){2}{\line(-1,-1){1}}

\multiput(0,4)(-1,-1){2}{\line(1,-1){1}}
\multiput(0,4)(1,-1){2}{\line(-1,-1){1}}
\multiput(4,4)(-1,-1){2}{\line(1,-1){1}}
\multiput(4,4)(1,-1){2}{\line(-1,-1){1}}

\multiput(1,3)(-1,-1){2}{\line(1,-1){1}}
\multiput(1,3)(1,-1){2}{\line(-1,-1){1}}
\multiput(3,3)(-1,-1){2}{\line(1,-1){1}}
\multiput(3,3)(1,-1){2}{\line(-1,-1){1}}

% e = 5, m = 6
\put(11,9.5){\makebox(0,0)[cb]{\(r=4\), \(c=5\)}}
\put(11,9){\circle*{0.2}}
\multiput(11,9)(-1,-1){2}{\line(1,-1){1}}
\multiput(11,9)(1,-1){2}{\line(-1,-1){1}}
\put(11,7){\circle*{0.2}}
% bottle neck
\put(11,7){\line(0,-1){1}}
%body
\put(11,6){\circle*{0.2}}
\multiput(11,6)(-1,-1){2}{\line(1,-1){1}}
\multiput(11,6)(1,-1){2}{\line(-1,-1){1}}
\put(11,4){\circle*{0.2}}
\multiput(11,4)(-1,-1){2}{\line(1,-1){1}}
\multiput(11,4)(1,-1){2}{\line(-1,-1){1}}
\put(11,2){\circle*{0.2}}
\multiput(11,2)(-1,-1){2}{\line(1,-1){1}}
\multiput(11,2)(1,-1){2}{\line(-1,-1){1}}
\put(11,0){\circle*{0.2}}

\multiput(10,5)(-1,-1){2}{\line(1,-1){1}}
\multiput(10,5)(1,-1){2}{\line(-1,-1){1}}
\multiput(12,5)(-1,-1){2}{\line(1,-1){1}}
\multiput(12,5)(1,-1){2}{\line(-1,-1){1}}

\multiput(9,4)(-1,-1){2}{\line(1,-1){1}}
\multiput(9,4)(1,-1){2}{\line(-1,-1){1}}
\multiput(13,4)(-1,-1){2}{\line(1,-1){1}}
\multiput(13,4)(1,-1){2}{\line(-1,-1){1}}

\multiput(10,3)(-1,-1){2}{\line(1,-1){1}}
\multiput(10,3)(1,-1){2}{\line(-1,-1){1}}
\multiput(12,3)(-1,-1){2}{\line(1,-1){1}}
\multiput(12,3)(1,-1){2}{\line(-1,-1){1}}

% e = 3, m = 4
\put(16,5.5){\makebox(0,0)[cb]{\(r=2\), \(c=3\)}}
\put(16,5){\circle*{0.2}}
\multiput(16,5)(-1,-1){2}{\line(1,-1){1}}
\multiput(16,5)(1,-1){2}{\line(-1,-1){1}}
\put(16,3){\circle*{0.2}}
% bottle neck
\put(16,3){\line(0,-1){1}}
%body
\put(16,2){\circle*{0.2}}
\multiput(16,2)(-1,-1){2}{\line(1,-1){1}}
\multiput(16,2)(1,-1){2}{\line(-1,-1){1}}
\put(16,0){\circle*{0.2}}

\end{picture}

\end{figure}

\noindent
In Figure
\ref{fig:BFGroups}
we display numerous examples of normal lattices of BF-groups \(G\)
with \textit{bicyclic factors} of the central series,
except the bottle neck \(\gamma_2(G)/\gamma_3(G)\).
They are located as vertices on the \textit{sporadic} part \(\mathcal{G}_0(p,r)\)
of coclass graphs \(\mathcal{G}(p,r)\) with \(r\ge 2\), outside of coclass trees,
\cite[Fig. 3.5, p. 439, Fig. 3.8, p. 448]{Ma2013}.

Here, the rectangle of trailing diamonds degenerates to a square with \(c=r+1\),
the upper central series is the reverse lower central series,
and thus the last lower central \(\gamma_{c}(G)\) is bicyclic,
whence the (generalized) parent \(\tilde\pi(G)=G/\gamma_{c}(G)\) is of lower coclass.
Such groups were called \textit{interface groups} to the coclass forest \(\mathcal{G}(p,r-1)\) in
\cite[Dfn. 3.3, p. 430]{Ma2013}.
They can also be realized as second \(p\)-class groups
\(M=\mathrm{Gal}(\mathrm{F}_p^2(K)\vert K)\)
by \textit{imaginary quadratic} fields \(K=\mathbb{Q}(\sqrt{d})\),
provided the coclass \(r\ge 2\) is \textit{even} (\textit{selection rule}),
and thus the nilpotency class \(c=r+1\) is \textit{odd}.

%--------------------------------------------------------------------------------

%--------------------------------------------------------------------------------

\begin{example}
\label{exm:BFGroups}
BF-groups \(M\) with \(p\in\lbrace 3,5,7\rbrace\).

\begin{itemize}
\item
From Table
\ref{tbl:Imag},
we get a rare example for 
\(p=3\), coclass \(\mathrm{cc}(M)=6\), class \(\mathrm{cl}(M)=7\),
realized by a single imaginary quadratic field in \(-10^6<d<0\), of discriminant\\
\(d=-423\,640\) with capitulation type (CT) F.12, and
a sporadic group \(M\) of order \(3^{13}=1\,594\,323\),
visualized by Figure
\ref{fig:BFGroups},
\(r=6\), \(c=7\).
\item
From
\cite[Tbl. 3, p. 497]{Ma2012a},
where the complete statistics for \(-10^6<d<0\) is listed,
we get numerous examples for
\(p=3\), coclass \(\mathrm{cc}(M)=4\), class \(\mathrm{cl}(M)=5\),
realized by a total of \(78\) imaginary quadratic fields, e.g.,\\
\(d=-27\,156\) with CT F.11,\\
\(d=-31\,908\) with CT F.12,\\
\(d=-67\,480\) with CT F.13,\\
\(d=-124\,363\) with CT F.7,\\
and by a single real quadratic field in \(0<d<10^7\), of discriminant\\
\(d=8\,321\,505\) with CT F.13,\\
with various non-isomorphic sporadic groups \(M\) of order \(3^9=19\,683\),
sharing isomophic normal lattices,
visualized by the same Figure
\ref{fig:BFGroups},
\(r=4\), \(c=5\).
\item
\(p=3\), coclass \(\mathrm{cc}(M)=2\), class \(\mathrm{cl}(M)=3\):\\
realized by a total of \(936\) imaginary quadratic fields in \(-10^6<d<0\), e.g.,\\
\(d=-4\,027\) with CT D.10,\\
\(d=-12\,131\) with CT D.5,\\
and by a total of \(140\) real quadratic fields in \(0<d<10^7\), e.g.,\\
\(d=422\,573\) with CT D.10,\\
\(d=631\,769\) with CT D.5,\\
two different sporadic Schur \(\sigma\)-groups,
\(\langle 243,5\rangle\) and \(\langle 243,7\rangle\),
sharing isomophic normal lattices,
visualized by the same Figure
\ref{fig:BFGroups},
\(r=2\), \(c=3\).
\item
\(p=5\), coclass \(\mathrm{cc}(M)=2\), class \(\mathrm{cl}(M)=3\): see
\cite[Tbl. 3.13, p. 450]{Ma2013}.
\item
\(p=7\), coclass \(\mathrm{cc}(M)=2\), class \(\mathrm{cl}(M)=3\): see
\cite[Tbl. 3.14, p. 450]{Ma2013}.
\end{itemize}

\end{example}

%\newpage
%--------------------------------------------------------------------------------

\subsection{Small bicyclic or cyclic factor groups}
\label{ss:SBCF}

\begin{figure}[ht]
\caption{Small BCF-groups \(M=\mathrm{Gal}(\mathrm{F}_p^2(K)\vert K)\) with bicyclic and cyclic factors.}
\label{fig:SmallBCFGroups}

% Small BCF Groups

\setlength{\unitlength}{0.5cm}
\begin{picture}(25,11)(-7,0)

% scale of orders
\put(-5.8,10.3){\makebox(0,0)[cb]{Order \(3^n\), \(p^n\)}}
\put(-5,9){\vector(0,1){1}}
\put(-5.2,9){\makebox(0,0)[rc]{\(19\,683\)}}
\put(-4.8,9){\makebox(0,0)[lc]{\(p^9\)}}
\put(-5.2,8){\makebox(0,0)[rc]{\(6\,561\)}}
\put(-4.8,8){\makebox(0,0)[lc]{\(p^8\)}}
\put(-5.2,7){\makebox(0,0)[rc]{\(2\,187\)}}
\put(-4.8,7){\makebox(0,0)[lc]{\(p^7\)}}
\put(-5.2,6){\makebox(0,0)[rc]{\(729\)}}
\put(-4.8,6){\makebox(0,0)[lc]{\(p^6\)}}
\put(-5.2,5){\makebox(0,0)[rc]{\(243\)}}
\put(-4.8,5){\makebox(0,0)[lc]{\(p^5\)}}
\put(-5.2,4){\makebox(0,0)[rc]{\(81\)}}
\put(-4.8,4){\makebox(0,0)[lc]{\(p^4\)}}
\put(-5.2,3){\makebox(0,0)[rc]{\(27\)}}
\put(-4.8,3){\makebox(0,0)[lc]{\(p^3\)}}
\put(-5.2,2){\makebox(0,0)[rc]{\(9\)}}
\put(-4.8,2){\makebox(0,0)[lc]{\(p^2\)}}
\put(-5.2,1){\makebox(0,0)[rc]{\(3\)}}
\put(-4.8,1){\makebox(0,0)[lc]{\(p\)}}
\put(-5.2,0){\makebox(0,0)[rc]{\(1\)}}
\multiput(-5.1,0)(0,1){10}{\line(1,0){0.2}}
\put(-5,0){\line(0,1){9}}

% e = 3, m = 8
\put(-2,9.5){\makebox(0,0)[cb]{\(r=2\), \(c=7\)}}
\put(-2,9){\circle*{0.2}}
\multiput(-2,9)(-1,-1){2}{\line(1,-1){1}}
\multiput(-2,9)(1,-1){2}{\line(-1,-1){1}}
\put(-2,7){\circle*{0.2}}
% bottle neck
\put(-2,7){\line(0,-1){1}}
% body
\put(-2,6){\circle*{0.2}}
\multiput(-2,6)(-1,-1){2}{\line(1,-1){1}}
\multiput(-2,6)(1,-1){2}{\line(-1,-1){1}}
\put(-2,4){\circle*{0.2}}

\put(-1,5){\circle{0.2}}
\multiput(-1,5)(-1,-1){2}{\line(1,-1){1}}
\multiput(-1,5)(1,-1){2}{\line(-1,-1){1}}
\put(-1,3){\circle*{0.2}}

\put(0,4){\circle{0.2}}
\multiput(0,4)(-1,-1){2}{\line(1,-1){1}}
\multiput(0,4)(1,-1){2}{\line(-1,-1){1}}
\put(0,2){\circle*{0.2}}

\put(1,3){\circle{0.2}}
\multiput(1,3)(-1,-1){2}{\line(1,-1){1}}
\multiput(1,3)(1,-1){2}{\line(-1,-1){1}}
\put(1,1){\circle*{0.2}}

\put(2,2){\circle{0.2}}
\multiput(2,2)(-1,-1){2}{\line(1,-1){1}}
\multiput(2,2)(1,-1){2}{\line(-1,-1){1}}
\put(2,0){\circle*{0.2}}

% e = 3, m = 7
\put(5,8.5){\makebox(0,0)[cb]{\(r=2\), \(c=6\)}}
\put(5,8){\circle*{0.2}}
\multiput(5,8)(-1,-1){2}{\line(1,-1){1}}
\multiput(5,8)(1,-1){2}{\line(-1,-1){1}}
\put(5,6){\circle*{0.2}}
% bottle neck
\put(5,6){\line(0,-1){1}}
% body
\put(5,5){\circle*{0.2}}
\multiput(5,5)(-1,-1){2}{\line(1,-1){1}}
\multiput(5,5)(1,-1){2}{\line(-1,-1){1}}
\put(5,3){\circle*{0.2}}

\put(6,4){\circle{0.2}}
\multiput(6,4)(-1,-1){2}{\line(1,-1){1}}
\multiput(6,4)(1,-1){2}{\line(-1,-1){1}}
\put(6,2){\circle*{0.2}}

\put(7,3){\circle{0.2}}
\multiput(7,3)(-1,-1){2}{\line(1,-1){1}}
\multiput(7,3)(1,-1){2}{\line(-1,-1){1}}
\put(7,1){\circle*{0.2}}

\put(8,2){\circle{0.2}}
\multiput(8,2)(-1,-1){2}{\line(1,-1){1}}
\multiput(8,2)(1,-1){2}{\line(-1,-1){1}}
\put(8,0){\circle*{0.2}}

% e = 3, m = 6
\put(11,7.5){\makebox(0,0)[cb]{\(r=2\), \(c=5\)}}
\put(11,7){\circle*{0.2}}
\multiput(11,7)(-1,-1){2}{\line(1,-1){1}}
\multiput(11,7)(1,-1){2}{\line(-1,-1){1}}
\put(11,5){\circle*{0.2}}
% bottle neck
\put(11,5){\line(0,-1){1}}
% body
\put(11,4){\circle*{0.2}}
\multiput(11,4)(-1,-1){2}{\line(1,-1){1}}
\multiput(11,4)(1,-1){2}{\line(-1,-1){1}}
\put(11,2){\circle*{0.2}}

\put(12,3){\circle{0.2}}
\multiput(12,3)(-1,-1){2}{\line(1,-1){1}}
\multiput(12,3)(1,-1){2}{\line(-1,-1){1}}
\put(12,1){\circle*{0.2}}

\put(13,2){\circle{0.2}}
\multiput(13,2)(-1,-1){2}{\line(1,-1){1}}
\multiput(13,2)(1,-1){2}{\line(-1,-1){1}}
\put(13,0){\circle*{0.2}}

% e = 3, m = 5
\put(16,6.5){\makebox(0,0)[cb]{\(r=2\), \(c=4\)}}
\put(16,6){\circle*{0.2}}
\multiput(16,6)(-1,-1){2}{\line(1,-1){1}}
\multiput(16,6)(1,-1){2}{\line(-1,-1){1}}
\put(16,4){\circle*{0.2}}
% bottle neck
\put(16,4){\line(0,-1){1}}
% body
\put(16,3){\circle*{0.2}}
\multiput(16,3)(-1,-1){2}{\line(1,-1){1}}
\multiput(16,3)(1,-1){2}{\line(-1,-1){1}}
\put(16,1){\circle*{0.2}}

\put(17,2){\circle{0.2}}
\multiput(17,2)(-1,-1){2}{\line(1,-1){1}}
\multiput(17,2)(1,-1){2}{\line(-1,-1){1}}
\put(17,0){\circle*{0.2}}

\end{picture}

\end{figure}
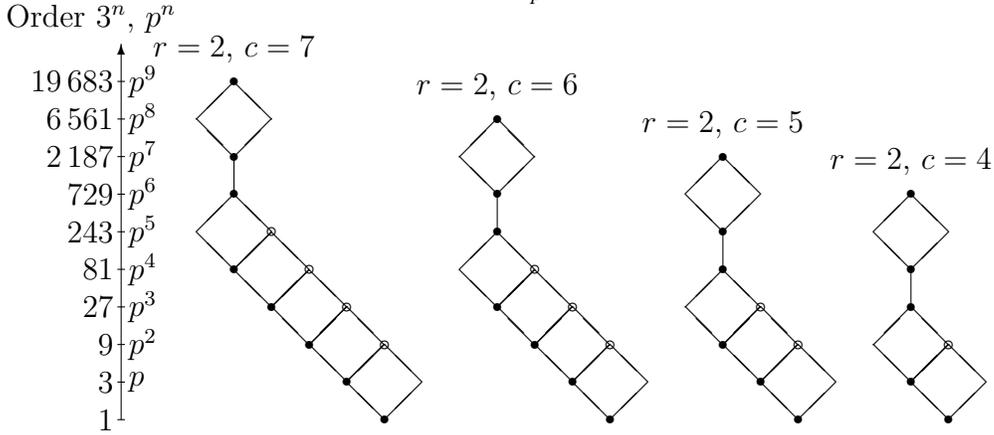

\noindent
Figure
\ref{fig:SmallBCFGroups}
shows many examples of normal lattices of \lq\lq small\rq\rq\ \(p\)-groups \(M\)
with \textit{bicyclic and cyclic factors} of the central series.
They are located on \textit{coclass trees} of the coclass graph \(\mathcal{G}(p,2)\)
\cite[Fig. 3.6--3.7, pp. 442--443]{Ma2013}.

%--------------------------------------------------------------------------------

%--------------------------------------------------------------------------------

\begin{example}
\label{exm:SmallBCFGroups}
Small BCF-groups \(M\) with \(p\in\lbrace 3,5,7\rbrace\), specified by capitulation types (CT).

\begin{itemize}
\item
\(p=3\), coclass \(\mathrm{cc}(M)=2\), class \(\mathrm{cl}(M)=7\):\\
a total of \(28\) imaginary quadratic fields \(K=\mathbb{Q}(\sqrt{d})\) in \(-10^6<d<0\), e.g.,\\
\(d=-262\,744\) with CT E.14\(\uparrow\), where \(\uparrow\) designates an excited state,\\
\(d=-268\,040\) with CT E.6\(\uparrow\),\\
\(d=-297\,079\) with CT E.9\(\uparrow\),\\
\(d=-370\,740\) with CT E.8\(\uparrow\),\\
\textit{branch} groups of depth \(1\), order \(3^9=19\,683\),
visualized by Figure
\ref{fig:SmallBCFGroups},
\(r=2\), \(c=7\).
\item
\(p=3\), coclass \(\mathrm{cc}(M)=2\), class \(\mathrm{cl}(M)=6\):\\
two real quadratic fields in \(0<d<10^7\), e. g.,\\
\(d=1\,001\,957\) with CT c.21\(\uparrow\), excited state,\\
\textit{mainline} groups of order \(3^8=6\,561\),
visualized by Figure
\ref{fig:SmallBCFGroups},
\(r=2\), \(c=6\).
\item
\(p=3\), coclass \(\mathrm{cc}(M)=2\), class \(\mathrm{cl}(M)=5\):\\
a total of \(383\) imaginary quadratic fields in \(-10^6<d<0\), e.g.,\\
\(d=-9\,748\) with CT E.9, ground state,\\
\(d=-15\,544\) with CT E.6,\\
\(d=-16\,627\) with CT E.14,\\
\(d=-34\,867\) with CT E.8,\\
and a total of \(21\) real quadratic fields in \(0<d<10^7\), e. g.,\\
\(d=342\,664\) with CT E.9, ground state,\\
\(d=3\,918\,837\) with CT E.14,\\
\(d=5\,264\,069\) with CT E.6,\\
\(d=6\,098\,360\) with CT E.8,\\
branch groups of depth \(1\), order \(3^7=2\,187\),
visualized by Figure
\ref{fig:SmallBCFGroups},
\(r=2\), \(c=5\).
\item
\(p=3\), coclass \(\mathrm{cc}(M)=2\), class \(\mathrm{cl}(M)=4\):\\
a total of \(53\) real quadratic fields in \(0<d<10^7\), e. g.,\\
\(d=534\,824\) with CT c.18, ground state,\\
\(d=540\,365\) with CT c.21,\\
mainline groups of order \(3^6=729\),
visualized by Figure
\ref{fig:SmallBCFGroups},
\(r=2\), \(c=4\).
\item
\(p=5\), coclass \(\mathrm{cc}(M)=2\), class \(\mathrm{cl}(M)=5\): see
\cite[Tbl. 3.13, p. 450]{Ma2013}.
\item
\(p=7\), coclass \(\mathrm{cc}(M)=2\), class \(\mathrm{cl}(M)=5\): see
\cite[Tbl. 3.14, p. 450]{Ma2013}.
\end{itemize}

\end{example}

%\newpage
%--------------------------------------------------------------------------------

\subsection{Large bicyclic or cyclic factor groups}
\label{ss:LBCF}

\noindent
Let \(K=\mathbb{Q}(\sqrt{d})\) be a quadratic number field with discriminant \(d\)
and denote by \(M=\mathrm{Gal}(\mathrm{F}_p^2(K)\vert K)\) the Galois group
of the second Hilbert \(p\)-class field \(\mathrm{F}_p^2(K)\) of \(K\),
that is, the maximal metabelian unramified \(p\)-extension of \(K\).
We recall that coclass and class of \(M\) are given by the equations
\(\mathrm{cc}(M)=r=e-1\) and \(\mathrm{cl}(M)=c=m-1\)
in terms of the CF-invariant \(e\) and the nilpotency index \(m\).
Due to our extensive computations for the papers
\cite{Ma2012a,Ma2013},
we are able to underpin the present theory of normal lattices
by numerical data concerning the \(2\,020\) imaginary and the \(2\,576\) real
quadratic fields with \(3\)-class group of type \((3,3)\)
and discriminants in the range \(-10^6<d<10^7\).

\begin{figure}[ht]
\caption{Large BCF-groups \(M=\mathrm{Gal}(\mathrm{F}_p^2(K)\vert K)\) with bicyclic and cyclic factors.}
\label{fig:BCFGroups}
% BCF Groups

\setlength{\unitlength}{0.5cm}
\begin{picture}(26,12)(-7,1)

% scale of orders
\put(-5.8,12.3){\makebox(0,0)[cb]{Order \(3^n\), \(p^n\)}}
\put(-5,11){\vector(0,1){1}}
\put(-5.2,11){\makebox(0,0)[rc]{\(177\,147\)}}
\put(-4.8,11){\makebox(0,0)[lc]{\(p^{11}\)}}
\put(-5.2,10){\makebox(0,0)[rc]{\(59\,049\)}}
\put(-4.8,10){\makebox(0,0)[lc]{\(p^{10}\)}}
\put(-5.2,9){\makebox(0,0)[rc]{\(19\,683\)}}
\put(-4.8,9){\makebox(0,0)[lc]{\(p^9\)}}
\put(-5.2,8){\makebox(0,0)[rc]{\(6\,561\)}}
\put(-4.8,8){\makebox(0,0)[lc]{\(p^8\)}}
\put(-5.2,7){\makebox(0,0)[rc]{\(2\,187\)}}
\put(-4.8,7){\makebox(0,0)[lc]{\(p^7\)}}
\put(-5.2,6){\makebox(0,0)[rc]{\(729\)}}
\put(-4.8,6){\makebox(0,0)[lc]{\(p^6\)}}
\put(-5.2,5){\makebox(0,0)[rc]{\(243\)}}
\put(-4.8,5){\makebox(0,0)[lc]{\(p^5\)}}
\put(-5.2,4){\makebox(0,0)[rc]{\(81\)}}
\put(-4.8,4){\makebox(0,0)[lc]{\(p^4\)}}
\put(-5.2,3){\makebox(0,0)[rc]{\(27\)}}
\put(-4.8,3){\makebox(0,0)[lc]{\(p^3\)}}
\put(-5.2,2){\makebox(0,0)[rc]{\(9\)}}
\put(-4.8,2){\makebox(0,0)[lc]{\(p^2\)}}
\put(-5.2,1){\makebox(0,0)[rc]{\(3\)}}
\put(-4.8,1){\makebox(0,0)[lc]{\(p\)}}
\put(-5.2,0){\makebox(0,0)[rc]{\(1\)}}
\multiput(-5.1,0)(0,1){12}{\line(1,0){0.2}}
\put(-5,0){\line(0,1){11}}

% e = 5, m = 8
\put(0,11.5){\makebox(0,0)[cb]{\(r=4\), \(c=7\)}}
\put(0,11){\circle*{0.2}}
\multiput(0,11)(-1,-1){2}{\line(1,-1){1}}
\multiput(0,11)(1,-1){2}{\line(-1,-1){1}}
\put(0,9){\circle*{0.2}}
% bottle neck
\put(0,9){\line(0,-1){1}}
% body
\put(0,8){\circle*{0.2}}
\multiput(0,8)(-1,-1){2}{\line(1,-1){1}}
\multiput(0,8)(1,-1){2}{\line(-1,-1){1}}
\put(0,6){\circle*{0.2}}
\multiput(0,6)(-1,-1){2}{\line(1,-1){1}}
\multiput(0,6)(1,-1){2}{\line(-1,-1){1}}
\put(0,4){\circle*{0.2}}
\multiput(0,4)(-1,-1){2}{\line(1,-1){1}}
\multiput(0,4)(1,-1){2}{\line(-1,-1){1}}
\put(0,2){\circle*{0.2}}

\put(1,7){\circle{0.2}}
\multiput(-1,7)(-1,-1){2}{\line(1,-1){1}}
\multiput(-1,7)(1,-1){2}{\line(-1,-1){1}}
\multiput(1,7)(-1,-1){2}{\line(1,-1){1}}
\multiput(1,7)(1,-1){2}{\line(-1,-1){1}}

\multiput(-2,6)(-1,-1){2}{\line(1,-1){1}}
\multiput(-2,6)(1,-1){2}{\line(-1,-1){1}}

\multiput(-1,5)(-1,-1){2}{\line(1,-1){1}}
\multiput(-1,5)(1,-1){2}{\line(-1,-1){1}}
\multiput(1,5)(-1,-1){2}{\line(1,-1){1}}
\multiput(1,5)(1,-1){2}{\line(-1,-1){1}}
\multiput(3,5)(-1,-1){2}{\line(1,-1){1}}
\multiput(3,5)(1,-1){2}{\line(-1,-1){1}}

\multiput(4,4)(-1,-1){2}{\line(1,-1){1}}
\multiput(4,4)(1,-1){2}{\line(-1,-1){1}}

\multiput(1,3)(-1,-1){2}{\line(1,-1){1}}
\multiput(1,3)(1,-1){2}{\line(-1,-1){1}}
\put(1,1){\circle*{0.2}}
\multiput(3,3)(-1,-1){2}{\line(1,-1){1}}
\multiput(3,3)(1,-1){2}{\line(-1,-1){1}}

\put(2,6){\circle{0.2}}
\multiput(2,6)(-1,-1){2}{\line(1,-1){1}}
\multiput(2,6)(1,-1){2}{\line(-1,-1){1}}
\put(2,4){\circle{0.2}}
\multiput(2,4)(-1,-1){2}{\line(1,-1){1}}
\multiput(2,4)(1,-1){2}{\line(-1,-1){1}}
\put(2,2){\circle{0.2}}
\multiput(2,2)(-1,-1){2}{\line(1,-1){1}}
\multiput(2,2)(1,-1){2}{\line(-1,-1){1}}
\put(2,0){\circle*{0.2}}

% e = 5, m = 7
\put(9,10.5){\makebox(0,0)[cb]{\(r=4\), \(c=6\)}}
\put(9,10){\circle*{0.2}}
\multiput(9,10)(-1,-1){2}{\line(1,-1){1}}
\multiput(9,10)(1,-1){2}{\line(-1,-1){1}}
\put(9,8){\circle*{0.2}}
% bottle neck
\put(9,8){\line(0,-1){1}}
% body
\multiput(8,6)(-1,-1){2}{\line(1,-1){1}}
\multiput(8,6)(1,-1){2}{\line(-1,-1){1}}

\multiput(7,5)(-1,-1){2}{\line(1,-1){1}}
\multiput(7,5)(1,-1){2}{\line(-1,-1){1}}
\multiput(11,5)(-1,-1){2}{\line(1,-1){1}}
\multiput(11,5)(1,-1){2}{\line(-1,-1){1}}

\multiput(8,4)(-1,-1){2}{\line(1,-1){1}}
\multiput(8,4)(1,-1){2}{\line(-1,-1){1}}
\multiput(12,4)(-1,-1){2}{\line(1,-1){1}}
\multiput(12,4)(1,-1){2}{\line(-1,-1){1}}

\multiput(11,3)(-1,-1){2}{\line(1,-1){1}}
\multiput(11,3)(1,-1){2}{\line(-1,-1){1}}

\put(9,7){\circle*{0.2}}
\multiput(9,7)(-1,-1){2}{\line(1,-1){1}}
\multiput(9,7)(1,-1){2}{\line(-1,-1){1}}
\put(9,5){\circle*{0.2}}
\multiput(9,5)(-1,-1){2}{\line(1,-1){1}}
\multiput(9,5)(1,-1){2}{\line(-1,-1){1}}
\put(9,3){\circle*{0.2}}
\multiput(9,3)(-1,-1){2}{\line(1,-1){1}}
\multiput(9,3)(1,-1){2}{\line(-1,-1){1}}
\put(9,1){\circle*{0.2}}

\put(10,6){\circle{0.2}}
\multiput(10,6)(-1,-1){2}{\line(1,-1){1}}
\multiput(10,6)(1,-1){2}{\line(-1,-1){1}}
\put(10,4){\circle{0.2}}
\multiput(10,4)(-1,-1){2}{\line(1,-1){1}}
\multiput(10,4)(1,-1){2}{\line(-1,-1){1}}
\put(10,2){\circle{0.2}}
\multiput(10,2)(-1,-1){2}{\line(1,-1){1}}
\multiput(10,2)(1,-1){2}{\line(-1,-1){1}}
\put(10,0){\circle*{0.2}}

% e = 4, m = 6
\put(16,8.5){\makebox(0,0)[cb]{\(r=3\), \(c=5\)}}
\put(16,8){\circle*{0.2}}
\multiput(16,8)(-1,-1){2}{\line(1,-1){1}}
\multiput(16,8)(1,-1){2}{\line(-1,-1){1}}
\put(16,6){\circle*{0.2}}
% bottle neck
\put(16,6){\line(0,-1){1}}
% body
\multiput(15,4)(-1,-1){2}{\line(1,-1){1}}
\multiput(15,4)(1,-1){2}{\line(-1,-1){1}}

\put(16,5){\circle*{0.2}}
\multiput(16,5)(-1,-1){2}{\line(1,-1){1}}
\multiput(16,5)(1,-1){2}{\line(-1,-1){1}}
\put(16,3){\circle*{0.2}}
\multiput(16,3)(-1,-1){2}{\line(1,-1){1}}
\multiput(16,3)(1,-1){2}{\line(-1,-1){1}}
\put(16,1){\circle*{0.2}}

\put(17,4){\circle{0.2}}
\multiput(17,4)(-1,-1){2}{\line(1,-1){1}}
\multiput(17,4)(1,-1){2}{\line(-1,-1){1}}
\put(17,2){\circle{0.2}}
\multiput(17,2)(-1,-1){2}{\line(1,-1){1}}
\multiput(17,2)(1,-1){2}{\line(-1,-1){1}}
\put(17,0){\circle*{0.2}}

\multiput(18,3)(-1,-1){2}{\line(1,-1){1}}
\multiput(18,3)(1,-1){2}{\line(-1,-1){1}}

\end{picture}

\end{figure}
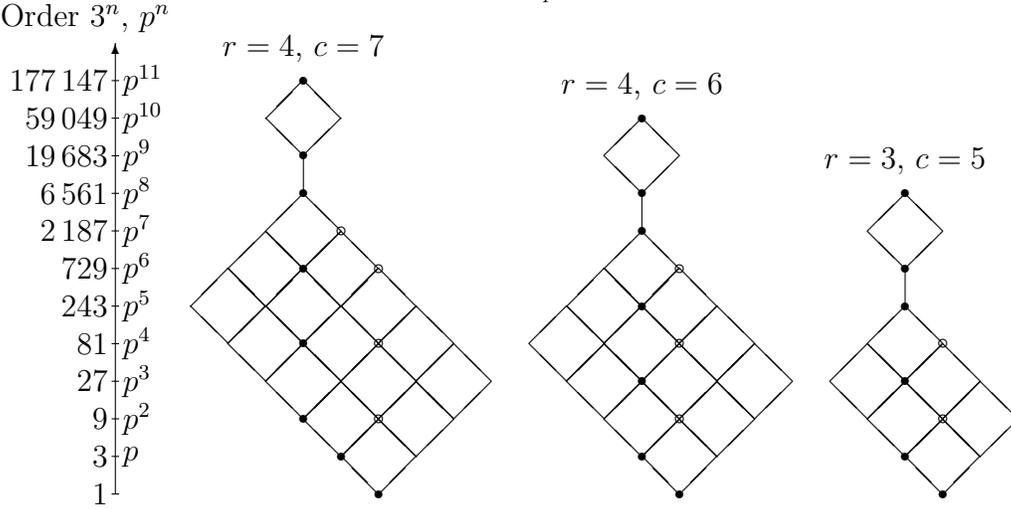

Figure
\ref{fig:BCFGroups}
shows several examples of normal lattices of \lq\lq large\rq\rq\ \(3\)-groups \(M\)
with \textit{bicyclic and cyclic factors} of the central series.
They are located on \textit{coclass trees} of coclass graphs \(\mathcal{G}(3,r)\) with \(r\ge 3\)
\cite[p. 189 ff]{Ne1989a}.

Here, the length of the rectangle of trailing diamonds is bigger than the width, \(c>r+1\),
the upper central series is \textit{different from} the lower central series,
and the last non-trivial lower central \(\gamma_{c}(M)\) is cyclic,
whence the parent \(\pi(M)=M/\gamma_{c}(M)\) is of the same coclass.
Such groups were called \textit{core groups} in
\cite{Ma2013}.
Concerning the capitulation type (CT) \(\varkappa(K)\) of \(K\),
which coincides with the transfer kernel type (TKT) \(\varkappa(M)\) of \(M\), see
\cite{Ma2012b,Ma2013}.
Different TKTs can give rise to isomorphic normal lattices.

%--------------------------------------------------------------------------------

%--------------------------------------------------------------------------------

\begin{example}
\label{exm:BCFGroups}
Large BCF-groups \(M\) of coclass \(3\le\mathrm{cc}(M)\le 6\), specified by capitulation types (CT).

\begin{itemize}
\item
Coclass \(\mathrm{cc}(M)=3\), class \(\mathrm{cl}(M)=5\):\\
realized by two real quadratic fields in \(0<d<10^7\), of discriminant\\
\(d=1\,535\,117\) with CT d.23,\\
\(d=2\,328\,721\) with CT d.19,\\
\textit{branch} groups of depth \(1\), order \(3^8=6\,561\),
visualized by Figure
\ref{fig:BCFGroups},
\(r=3\), \(c=5\).
\item
Coclass \(\mathrm{cc}(M)=3\), class \(\mathrm{cl}(M)=7\):\\
realized by a real quadratic field of minimal discriminant
\(d=27\,970\,737\) with CT d.19
(not drawn in Figure
\ref{fig:BCFGroups}).
\item
Coclass \(\mathrm{cc}(M)=3\), class \(\mathrm{cl}(M)=9\):\\
realized by a real quadratic field of minimal discriminant
\(d=131\,279\,821\) with CT d.25
(not drawn in Figure
\ref{fig:BCFGroups}).
\item
Coclass \(\mathrm{cc}(M)=4\), class \(\mathrm{cl}(M)=6\):\\
realized by a single real quadratic field in \(0<d<10^7\), of discriminant\\
\(d=8\,491\,713\) with CT d.25\({}^\ast\),\\
\textit{mainline} group of order \(3^{10}=59\,049\),
visualized by Figure
\ref{fig:BCFGroups},
\(r=4\), \(c=6\).
\item
Coclass \(\mathrm{cc}(M)=4\), class \(\mathrm{cl}(M)=7\):\\
realized by a total of \(14\) imaginary quadratic fields in \(-10^6<d<0\), e. g.,\\
\(d=-159\,208\) with CT F.13,\\
\(d=-249\,371\) with CT F.12,\\
\(d=-469\,787\) with CT F.11,\\
\(d=-469\,816\) with CT F.7,\\
and by a single real quadratic field in \(0<d<10^7\), of discriminant\\
\(d=8\,127\,208\) with CT F.13,\\
\textit{branch} groups of depth \(1\), order \(3^{11}=177\,147\),
visualized by Figure
\ref{fig:BCFGroups},
\(r=4\), \(c=7\).
\item
Coclass \(\mathrm{cc}(M)=4\), class \(\mathrm{cl}(M)=9\):\\
realized by an imaginary quadratic field of absolutely minimal discriminant
\(d=-1\,699\,711\) with CT F.12
(not drawn in Figure
\ref{fig:BCFGroups}).
\item
Coclass \(\mathrm{cc}(M)=4\), class \(\mathrm{cl}(M)=11\):\\
realized by an imaginary quadratic field of absolutely minimal discriminant
\(d=-21\,548\,328\) with CT F.12
(not drawn in Figure
\ref{fig:BCFGroups}).
\item
Coclass \(\mathrm{cc}(M)=6\), class \(\mathrm{cl}(M)=9\):\\
realized by an imaginary quadratic field of absolutely minimal discriminant
\(d=-3\,820\,387\) with CT F.7
(not drawn in Figure
\ref{fig:BCFGroups}).
\item
Coclass \(\mathrm{cc}(M)=6\), class \(\mathrm{cl}(M)=11\):\\
realized by an imaginary quadratic field of absolutely minimal discriminant
\(d=-72\,411\,587\) with CT F.12
(not drawn in Figure
\ref{fig:BCFGroups}).
\end{itemize}

\end{example}

%\newpage
%--------------------------------------------------------------------------------

\subsection{The current coclass record}
\label{ss:RecordCoclass}

\begin{figure}[ht]
\caption{Biggest BCF-normal lattice currently realized arithmetically.}
\label{fig:BiggestLattice}

% BF Groups

\setlength{\unitlength}{0.5cm}
\begin{picture}(24,20)(-10,1)

% scale of orders
\put(-5.8,20.3){\makebox(0,0)[cb]{Order \(3^n\), \(p^n\)}}
\put(-5,19){\vector(0,1){1}}
\put(-5.2,19){\makebox(0,0)[rc]{\(1\,162\,261\,467\)}}
\put(-4.8,19){\makebox(0,0)[lc]{\(p^{19}\)}}
\put(-5.2,18){\makebox(0,0)[rc]{\(387\,420\,489\)}}
\put(-4.8,18){\makebox(0,0)[lc]{\(p^{18}\)}}
\put(-5.2,17){\makebox(0,0)[rc]{\(129\,140\,163\)}}
\put(-4.8,17){\makebox(0,0)[lc]{\(p^{17}\)}}
\put(-5.2,16){\makebox(0,0)[rc]{\(43\,046\,721\)}}
\put(-4.8,16){\makebox(0,0)[lc]{\(p^{16}\)}}
\put(-5.2,15){\makebox(0,0)[rc]{\(14\,348\,907\)}}
\put(-4.8,15){\makebox(0,0)[lc]{\(p^{15}\)}}
\put(-5.2,14){\makebox(0,0)[rc]{\(4\,782\,969\)}}
\put(-4.8,14){\makebox(0,0)[lc]{\(p^{14}\)}}
\put(-5.2,13){\makebox(0,0)[rc]{\(1\,594\,323\)}}
\put(-4.8,13){\makebox(0,0)[lc]{\(p^{13}\)}}
\put(-5.2,12){\makebox(0,0)[rc]{\(531\,441\)}}
\put(-4.8,12){\makebox(0,0)[lc]{\(p^{12}\)}}
\put(-5.2,11){\makebox(0,0)[rc]{\(177\,147\)}}
\put(-4.8,11){\makebox(0,0)[lc]{\(p^{11}\)}}
\put(-5.2,10){\makebox(0,0)[rc]{\(59\,049\)}}
\put(-4.8,10){\makebox(0,0)[lc]{\(p^{10}\)}}
\put(-5.2,9){\makebox(0,0)[rc]{\(19\,683\)}}
\put(-4.8,9){\makebox(0,0)[lc]{\(p^9\)}}
\put(-5.2,8){\makebox(0,0)[rc]{\(6\,561\)}}
\put(-4.8,8){\makebox(0,0)[lc]{\(p^8\)}}
\put(-5.2,7){\makebox(0,0)[rc]{\(2\,187\)}}
\put(-4.8,7){\makebox(0,0)[lc]{\(p^7\)}}
\put(-5.2,6){\makebox(0,0)[rc]{\(729\)}}
\put(-4.8,6){\makebox(0,0)[lc]{\(p^6\)}}
\put(-5.2,5){\makebox(0,0)[rc]{\(243\)}}
\put(-4.8,5){\makebox(0,0)[lc]{\(p^5\)}}
\put(-5.2,4){\makebox(0,0)[rc]{\(81\)}}
\put(-4.8,4){\makebox(0,0)[lc]{\(p^4\)}}
\put(-5.2,3){\makebox(0,0)[rc]{\(27\)}}
\put(-4.8,3){\makebox(0,0)[lc]{\(p^3\)}}
\put(-5.2,2){\makebox(0,0)[rc]{\(9\)}}
\put(-4.8,2){\makebox(0,0)[lc]{\(p^2\)}}
\put(-5.2,1){\makebox(0,0)[rc]{\(3\)}}
\put(-4.8,1){\makebox(0,0)[lc]{\(p\)}}
\put(-5.2,0){\makebox(0,0)[rc]{\(1\)}}
\multiput(-5.1,0)(0,1){20}{\line(1,0){0.2}}
\put(-5,0){\line(0,1){19}}

% e = 9, m = 12
\put(4,19.5){\makebox(0,0)[cb]{\(r=8\), \(c=11\)}}
\put(4,19){\circle*{0.2}}
\multiput(4,19)(-1,-1){2}{\line(1,-1){1}}
\multiput(4,19)(1,-1){2}{\line(-1,-1){1}}
\put(4,17){\circle*{0.2}}
% bottle neck
\put(4,17){\line(0,-1){1}}
% body
\put(4,16){\circle*{0.2}}
\multiput(4,16)(-1,-1){2}{\line(1,-1){1}}
\multiput(4,16)(1,-1){2}{\line(-1,-1){1}}
\put(4,14){\circle*{0.2}}
\multiput(4,14)(-1,-1){2}{\line(1,-1){1}}
\multiput(4,14)(1,-1){2}{\line(-1,-1){1}}
\put(4,12){\circle*{0.2}}
\multiput(4,12)(-1,-1){2}{\line(1,-1){1}}
\multiput(4,12)(1,-1){2}{\line(-1,-1){1}}
\put(4,10){\circle*{0.2}}
\multiput(4,10)(-1,-1){2}{\line(1,-1){1}}
\multiput(4,10)(1,-1){2}{\line(-1,-1){1}}
\put(4,8){\circle*{0.2}}
\multiput(4,8)(-1,-1){2}{\line(1,-1){1}}
\multiput(4,8)(1,-1){2}{\line(-1,-1){1}}
\put(4,6){\circle*{0.2}}
\multiput(4,6)(-1,-1){2}{\line(1,-1){1}}
\multiput(4,6)(1,-1){2}{\line(-1,-1){1}}
\put(4,4){\circle*{0.2}}
\multiput(4,4)(-1,-1){2}{\line(1,-1){1}}
\multiput(4,4)(1,-1){2}{\line(-1,-1){1}}
\put(4,2){\circle*{0.2}}

\put(6,14){\circle{0.2}}
\multiput(6,14)(-1,-1){2}{\line(1,-1){1}}
\multiput(6,14)(1,-1){2}{\line(-1,-1){1}}
\put(6,12){\circle{0.2}}
\multiput(6,12)(-1,-1){2}{\line(1,-1){1}}
\multiput(6,12)(1,-1){2}{\line(-1,-1){1}}
\put(6,10){\circle{0.2}}
\multiput(6,10)(-1,-1){2}{\line(1,-1){1}}
\multiput(6,10)(1,-1){2}{\line(-1,-1){1}}
\put(6,8){\circle{0.2}}
\multiput(6,8)(-1,-1){2}{\line(1,-1){1}}
\multiput(6,8)(1,-1){2}{\line(-1,-1){1}}
\put(6,6){\circle{0.2}}
\multiput(6,6)(-1,-1){2}{\line(1,-1){1}}
\multiput(6,6)(1,-1){2}{\line(-1,-1){1}}
\put(6,4){\circle{0.2}}
\multiput(6,4)(-1,-1){2}{\line(1,-1){1}}
\multiput(6,4)(1,-1){2}{\line(-1,-1){1}}
\put(6,2){\circle{0.2}}
\multiput(6,2)(-1,-1){2}{\line(1,-1){1}}
\multiput(6,2)(1,-1){2}{\line(-1,-1){1}}
\put(6,0){\circle*{0.2}}

\multiput(3,15)(-1,-1){2}{\line(1,-1){1}}
\multiput(3,15)(1,-1){2}{\line(-1,-1){1}}
\put(5,15){\circle{0.2}}
\multiput(5,15)(-1,-1){2}{\line(1,-1){1}}
\multiput(5,15)(1,-1){2}{\line(-1,-1){1}}

\multiput(2,14)(-1,-1){2}{\line(1,-1){1}}
\multiput(2,14)(1,-1){2}{\line(-1,-1){1}}

\multiput(1,13)(-1,-1){2}{\line(1,-1){1}}
\multiput(1,13)(1,-1){2}{\line(-1,-1){1}}
\multiput(3,13)(-1,-1){2}{\line(1,-1){1}}
\multiput(3,13)(1,-1){2}{\line(-1,-1){1}}
\multiput(5,13)(-1,-1){2}{\line(1,-1){1}}
\multiput(5,13)(1,-1){2}{\line(-1,-1){1}}
\multiput(7,13)(-1,-1){2}{\line(1,-1){1}}
\multiput(7,13)(1,-1){2}{\line(-1,-1){1}}

\multiput(0,12)(-1,-1){2}{\line(1,-1){1}}
\multiput(0,12)(1,-1){2}{\line(-1,-1){1}}
\multiput(2,12)(-1,-1){2}{\line(1,-1){1}}
\multiput(2,12)(1,-1){2}{\line(-1,-1){1}}
\multiput(8,12)(-1,-1){2}{\line(1,-1){1}}
\multiput(8,12)(1,-1){2}{\line(-1,-1){1}}

\multiput(-1,11)(-1,-1){2}{\line(1,-1){1}}
\multiput(-1,11)(1,-1){2}{\line(-1,-1){1}}
\multiput(1,11)(-1,-1){2}{\line(1,-1){1}}
\multiput(1,11)(1,-1){2}{\line(-1,-1){1}}
\multiput(3,11)(-1,-1){2}{\line(1,-1){1}}
\multiput(3,11)(1,-1){2}{\line(-1,-1){1}}
\multiput(5,11)(-1,-1){2}{\line(1,-1){1}}
\multiput(5,11)(1,-1){2}{\line(-1,-1){1}}
\multiput(7,11)(-1,-1){2}{\line(1,-1){1}}
\multiput(7,11)(1,-1){2}{\line(-1,-1){1}}
\multiput(9,11)(-1,-1){2}{\line(1,-1){1}}
\multiput(9,11)(1,-1){2}{\line(-1,-1){1}}

\multiput(-2,10)(-1,-1){2}{\line(1,-1){1}}
\multiput(-2,10)(1,-1){2}{\line(-1,-1){1}}
\multiput(0,10)(-1,-1){2}{\line(1,-1){1}}
\multiput(0,10)(1,-1){2}{\line(-1,-1){1}}
\multiput(2,10)(-1,-1){2}{\line(1,-1){1}}
\multiput(2,10)(1,-1){2}{\line(-1,-1){1}}
\multiput(8,10)(-1,-1){2}{\line(1,-1){1}}
\multiput(8,10)(1,-1){2}{\line(-1,-1){1}}
\multiput(10,10)(-1,-1){2}{\line(1,-1){1}}
\multiput(10,10)(1,-1){2}{\line(-1,-1){1}}

\multiput(-1,9)(-1,-1){2}{\line(1,-1){1}}
\multiput(-1,9)(1,-1){2}{\line(-1,-1){1}}
\multiput(1,9)(-1,-1){2}{\line(1,-1){1}}
\multiput(1,9)(1,-1){2}{\line(-1,-1){1}}
\multiput(3,9)(-1,-1){2}{\line(1,-1){1}}
\multiput(3,9)(1,-1){2}{\line(-1,-1){1}}
\multiput(5,9)(-1,-1){2}{\line(1,-1){1}}
\multiput(5,9)(1,-1){2}{\line(-1,-1){1}}
\multiput(7,9)(-1,-1){2}{\line(1,-1){1}}
\multiput(7,9)(1,-1){2}{\line(-1,-1){1}}
\multiput(9,9)(-1,-1){2}{\line(1,-1){1}}
\multiput(9,9)(1,-1){2}{\line(-1,-1){1}}
\multiput(11,9)(-1,-1){2}{\line(1,-1){1}}
\multiput(11,9)(1,-1){2}{\line(-1,-1){1}}

\multiput(0,8)(-1,-1){2}{\line(1,-1){1}}
\multiput(0,8)(1,-1){2}{\line(-1,-1){1}}
\multiput(2,8)(-1,-1){2}{\line(1,-1){1}}
\multiput(2,8)(1,-1){2}{\line(-1,-1){1}}
\multiput(8,8)(-1,-1){2}{\line(1,-1){1}}
\multiput(8,8)(1,-1){2}{\line(-1,-1){1}}
\multiput(10,8)(-1,-1){2}{\line(1,-1){1}}
\multiput(10,8)(1,-1){2}{\line(-1,-1){1}}
\multiput(12,8)(-1,-1){2}{\line(1,-1){1}}
\multiput(12,8)(1,-1){2}{\line(-1,-1){1}}

\multiput(1,7)(-1,-1){2}{\line(1,-1){1}}
\multiput(1,7)(1,-1){2}{\line(-1,-1){1}}
\multiput(3,7)(-1,-1){2}{\line(1,-1){1}}
\multiput(3,7)(1,-1){2}{\line(-1,-1){1}}
\multiput(7,7)(-1,-1){2}{\line(1,-1){1}}
\multiput(7,7)(1,-1){2}{\line(-1,-1){1}}
\multiput(9,7)(-1,-1){2}{\line(1,-1){1}}
\multiput(9,7)(1,-1){2}{\line(-1,-1){1}}
\multiput(11,7)(-1,-1){2}{\line(1,-1){1}}
\multiput(11,7)(1,-1){2}{\line(-1,-1){1}}

\multiput(2,6)(-1,-1){2}{\line(1,-1){1}}
\multiput(2,6)(1,-1){2}{\line(-1,-1){1}}
\multiput(8,6)(-1,-1){2}{\line(1,-1){1}}
\multiput(8,6)(1,-1){2}{\line(-1,-1){1}}
\multiput(10,6)(-1,-1){2}{\line(1,-1){1}}
\multiput(10,6)(1,-1){2}{\line(-1,-1){1}}

\multiput(3,5)(-1,-1){2}{\line(1,-1){1}}
\multiput(3,5)(1,-1){2}{\line(-1,-1){1}}
\multiput(7,5)(-1,-1){2}{\line(1,-1){1}}
\multiput(7,5)(1,-1){2}{\line(-1,-1){1}}
\multiput(9,5)(-1,-1){2}{\line(1,-1){1}}
\multiput(9,5)(1,-1){2}{\line(-1,-1){1}}

\multiput(8,4)(-1,-1){2}{\line(1,-1){1}}
\multiput(8,4)(1,-1){2}{\line(-1,-1){1}}

\multiput(5,3)(-1,-1){2}{\line(1,-1){1}}
\multiput(5,3)(1,-1){2}{\line(-1,-1){1}}
\put(5,1){\circle*{0.2}}
\multiput(7,3)(-1,-1){2}{\line(1,-1){1}}
\multiput(7,3)(1,-1){2}{\line(-1,-1){1}}

\end{picture}

\end{figure}
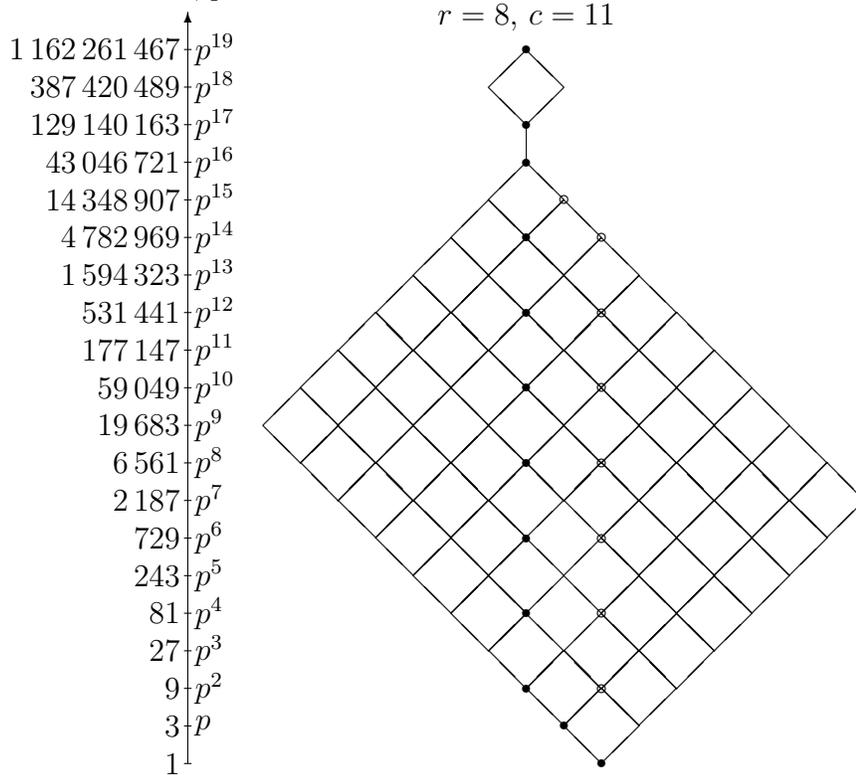

\noindent
At the moment,
the experimental state of the art is
a maximum of the coclass \(\mathrm{cc}(M)=8\),
realized by the unique imaginary quadratic field with discriminant
\(d=-99\,888\,340\)
from Table
\ref{tbl:Imag}.
Since the nilpotency class \(\mathrm{cl}(M)=11>9\) is not minimal,
the second \(3\)-class group \(M\) is not an interface BF-group,
and the normal lattice of this monster group, with order
\(\#M=3^{19}=1\,162\,261\,467\) and CT F.13,
has trailing diamonds in a proper rectangle, rather than a square.
This core BCF-group lattice is drawn in Figure
\ref{fig:BiggestLattice}. 

%--------------------------------------------------------------------------------

%\newpage
%--------------------------------------------------------------------------------

\subsection{Final Remarks}
\label{ss:FinalRemarks}

\begin{itemize}

\item
Among the \(2\,020\) imaginary quadratic fields
with \(3\)-class group of type \((3,3)\)
and discriminant in the range \(-10^6<d<0\),
the dominating part of \(1\,440\), that is \(71.29\,\%\),
has second \(3\)-class group
with minimal defect of commutativity \(k(M)=0\).
The remaining \(28.71\,\%\) have \(k(M)=1\)
and TKTs G.16, G.19 and H.4.

\item
Among the \(2\,576\) real quadratic fields
with \(3\)-class group of type \((3,3)\)
and discriminant in the range \(0<d<10^7\),
a modest part of \(273\), i. e. \(10.6\,\%\),
has a second \(3\)-class group of coclass at least \(2\).
A dominating part of \(222\) among these \(273\) second \(3\)-class groups,
that is \(81.3\,\%\),
has minimal defect of commutativity \(k(M)=0\),
whereas \(18.7\,\%\) have \(k(M)=1\)
and TKTs b.10, G.16, G.19 and H.4.

\item
It should be pointed out that
the power-commutator presentations which we used for proving Theorem
\ref{thm:NormalLattice}
and its Corollaries are rudimentary,
since in fact they consist of commutator relations only.
Thus they define an isoclinism family of \(p\)-groups of fixed order,
rather than a single isomorphism class of \(p\)-groups.

On the other hand,
experience shows that the transfer kernel type (TKT) of a \(p\)-group
mainly depends on the power relations.
This explains why different TKTs frequently give rise to equal normal lattices.

\end{itemize}

%\newpage
%--------------------------------------------------------------------------------

\section{Historical Overview}
\label{s:History}

\noindent
Motivated by progress in the classification
of \(p\)-groups of maximal class by Blackburn in \(1958\)
\cite{Bl1958},
and \(3\)-groups of second maximal class by Ascione et al. in \(1977\)
\cite{AHL1977},
Leedham-Green and Newman
\cite{LgNm1980}
proposed in \(1980\) to use the coclass \(\mathrm{cc}(G)\) as the primary invariant
for the classification of finite \(p\)-groups \(G\). 
The importance of the coclass as an invariant of finite \(3\)-groups
began to attract the attention and vigilance of scientists in \(1989\),
due to the Ph.D. thesis
\cite{Ne1989a,Ne1989b}
of Brigitte Nebelung,
since a periodicity with respect to the nilpotency class \(\mathrm{cl}(M)\)
and a co-periodicity with respect to the coclass \(\mathrm{cc}(M)\)
admitted to establish a finite parametrized set of power-commutator-presentations
for all metabelian \(3\)-groups \(M\) with
elementary bicyclic commutator quotient \(M/M^\prime\simeq (3,3)\).

For real quadratic fields \(K\),
the first progress beyond the prototype \(d=32\,009\) with coclass \(\mathrm{cc}(M)=1\)
by Heider and Schmithals \cite{HeSm1982}
was achieved on 28 January \(2006\) by our discovery of the discriminant \(d=214\,712\)
with coclass \(\mathrm{cc}(M)=2\).
On 09 November \(2009\), we succeeded in finding the 
minimal discriminant \(d=710\,652\)
with coclass \(\mathrm{cc}(M)=3\).
Until the end of \(2009\) we used the
cumbersome algorithm of Scholz and Taussky \cite{SoTa1934}
to compute the capitulation type \(\varkappa(K)\).
However, at the beginning of \(2010\)
we developed a marvellous technique
\cite{Ma2014a}
which allows to conclude the transfer kernel type \(\varkappa(K)\)
from the transfer target type \(\tau(K)\).
This yielded \(d=8\,127\,208\)
with \(\mathrm{cc}(M)=4\) on 18 March \(2010\).

%\newpage
%--------------------------------------------------------------------------------

\section{Outlook}
\label{s:Outlook}

\noindent
In order to extend our Tables
\ref{tbl:Imag} and \ref{tbl:Real}
to the minimal discriminants 
\(\lvert d\rvert\) for coclass \(\mathrm{cc}(M)=10\)
over imaginary quadratic fields, and
\(d\) for coclass \(\mathrm{cc}(M)=6\)
over real quadratic fields,
a gapless list of IPADs
\(\tau_1(K)=\left(\mathrm{ATI}(\mathrm{Cl}_3(E_i))\right)_{1\le i\le 4}\)
over quadratic fields \(K\)
with an estimated range \(-10^{11}<d<10^{12}\) of discriminants
would be required.
Since Michael Raymond Bush
\cite{BBH2017,BBH2021}
computed the most extensive range \(-10^{8}<d<10^{9}\) up to now
on a cluster of highly parallel super computers
with a total CPU time of several months,
the extension seems out of reach currently.

%\newpage
%--------------------------------------------------------------------------------

\section{Data availability statement}
\label{s:Data}

\noindent
Experimental results communicated in this article and
the source code of Magma scripts
\cite{BCP1997,BCFS2025,MAGMA2025}
used for the computations
may be requested from the second author by email. 

%\newpage
%--------------------------------------------------------------------------------

\section{Acknowledgements}
\label{s:Acknowledgements}

\noindent
The second author acknowledges that his research was supported by
%The author acknowledges that his research was supported by
the Austrian Science Fund (FWF): projects J0497-PHY, P26008-N25,
and by the Research Executive Agency of the European Union (EUREA):
project Horizon Europe 2021--2027.

%\newpage
%--------------------------------------------------------------------------------

%--------------------------------------------------------------------------------


\begin{thebibliography}{XX}
%
\bibitem{AoMa2024}
S. Aouissi and D. C. Mayer,
\textit{A group theoretic approach to cyclic cubic fields},
Mathematics
\textbf{12}
(2024),
no. 126,
1--52,
DOI 10.3390/math12010126.
%
\bibitem{AoMa2025}
S. Aouissi and D. C. Mayer,
\textit{The capitulation problem in certain pure cubic fields},
accepted by Annales Math\'ematiques du Qu\'ebec,
2025.
%
\bibitem{AHL1977}
J. A. Ascione, G. Havas, and C. R. Leedham-Green,
\textit{A computer aided classification of certain groups of prime power order},
Bull. Austral. Math. Soc.
\textbf{17}
(1977),
257--274, Corrigendum 317--319, Microfiche Supplement p. 320.
%
\bibitem{BEO2002}
H. U. Besche, B. Eick, and E. A. O'Brien,
\textit{A millennium project: constructing small groups},
Int. J. Algebra Comput.
\textbf{12}
(2002),
623-644,
DOI 10.1142/s0218196702001115.
%
\bibitem{BEO2005}
H. U. Besche, B. Eick, and E. A. O'Brien,
\textit{The SmallGroups Library --- a Library of Groups of Small Order},
2005,
an accepted and refereed GAP package, available also in Magma.
%
\bibitem{Bl1958}
N. Blackburn,
\textit{On a special class of \(p\)-groups},
Acta Math.
\textbf{100}
(1958),
45--92.
%
\bibitem{BCP1997}
W. Bosma, J. Cannon, and C. Playoust,
\textit{The Magma algebra system. I. The user language},
J. Symbolic Comput.
\textbf{24}
(1997),
235--265.
%
\bibitem{BCFS2025}
W. Bosma, J. J. Cannon, C. Fieker, and A. Steel (eds.),
\textit{Handbook of Magma functions},
Edition 2.28,
Sydney,
2025.
%
\bibitem{BBH2017}
N. Boston, M. R. Bush and F. Hajir,
\textit{Heuristics for \(p\)-class towers of imaginary quadratic fields},
Math. Ann.
\textbf{368}
(2017),
no. 1,
633--669,
DOI 10.1007/s00208-016-1449-3.
%
\bibitem{BBH2021}
N. Boston, M. R. Bush and F. Hajir,
\textit{Heuristics for \(p\)-class towers of real quadratic fields},
J. Inst. Math. Jussieu
\textbf{20}
(2021),
no. 4,
1429--1452,
DOI: 10.1017/S1474748019000641.
%
\bibitem{Fi2001}
C. Fieker,
\textit{Computing class fields via the Artin map},
Math. Comp.
\textbf{70}
(2001),
no. 235,
1293--1303.
%
\bibitem{GNO2006}
G. Gamble, W. Nickel, and E. A. O'Brien,
\textit{ANU p-Quotient --- p-Quotient and p-Group Generation Algorithms},
2006,
an accepted GAP package,
available also in Magma.
%
\bibitem{GAP2025}
GAP Developer Group,
GAP -- Groups, Algorithms, and Programming, Version \texttt{4.14.0},
2024,
available from
\texttt{http://www.gap-system.org}.
%
\bibitem{HeSm1982}
F.-P. Heider und B. Schmithals,
\textit{Zur Kapitulation der Idealklassen in unverzweigten primzyklischen Erweiterungen},
J. Reine Angew. Math.
\textbf{336}
(1982),
1--25.
%
\bibitem{HEO2005}
D. F. Holt, B. Eick, and E. A. O'Brien,
\textit{Handbook of computational group theory},
Discrete mathematics and its applications,
Chapman and Hall/CRC Press,
2005.
%
\bibitem{Hu1967}
B. Huppert,
Endliche Gruppen, I,
Springer-Verlag, Berlin,
1967.
%
\bibitem{Jm1980}
R. James,
\textit{The groups of order \(p^6\) (\(p\) an odd prime).},
Math. Comp
\textbf{34}
(1980),
no. 150,
613--637.
%
\bibitem{LgNm1980}
C. R. Leedham-Green and M. F. Newman,
\textit{Space groups and groups of prime power order I},
Arch. Math.
\textbf{35}
(1980),
193--203.
%
\bibitem{MAGMA2025}
Magma Developer Group,
\textit{Magma Computational Algebra System},
Version \texttt{2.28-27},
Sydney,
2025,
available from
\texttt{http://magma.maths.usyd.edu.au}.
%
\bibitem{MAGMA6561}
Magma Developer Group,
\textit{Magma, Data for groups of order \(3^8\)},
\texttt{data3to8.tar.gz},
Sydney,
2012,
available from
\texttt{http://magma.maths.usyd.edu.au}.
%
\bibitem{Ma1992}
D. C. Mayer, 
\textit{Multiplicities of dihedral discriminants}, 
Math. Comp.
\textbf{58}
(1992),
no. 198,
831--847,
supplements section S55--S58.
%
\bibitem{Ma2012a}
D. C. Mayer,
\textit{The second \(p\)-class group of a number field},
Int. J. Number Theory
\textbf{8}
(2012),
no. 2,
471--505,
DOI 10.1142/S179304211250025X.
%
\bibitem{Ma2012b}
D. C. Mayer,
\textit{Transfers of metabelian \(p\)-groups},
Monatsh. Math.
\textbf{166}
(2012),
no. 3--4,
467--495,
DOI 10.1007/s00605-010-0277-x.
%
\bibitem{Ma2013}
D. C. Mayer,
\textit{The distribution of second \(p\)-class groups on coclass graphs},
J. Th\'eor. Nombres Bordeaux
\textbf{25}
(2013),
no. 2,
401--456,
DOI 10.5802/jtnb842.
%
\bibitem{Ma2014a}
D. C. Mayer,
\textit{Principalization algorithm via class group structure},
J. Th\'eor. Nombres Bordeaux
\textbf{26}
(2014),
no. 2,
415--464,
DOI 10.5802/jtnb.874.
%
\bibitem{Ma2014b}
D. C. Mayer,
\textit{Quadratic \(p\)-ring spaces for counting dihedral fields},
Int. J. Number Theory
\textbf{10}
(2014),
no. 8,
2205--2242,
DOI 10.1142/S1793042114500754.
%
\bibitem{Ma2015a}
D. C. Mayer,
\textit{Periodic bifurcations in descendant trees of finite \(p\)-groups},
Adv. Pure Math.
\textbf{5}
(2015),
no. 4,
162--195,
DOI 10.4236/apm.2015.54020.
%
\bibitem{Ma2015b}
D. C. Mayer,
\textit{Index-\(p\) abelianization data of \(p\)-class tower groups},
Adv. Pure Math.
\textbf{5}
(2015)
No. 5,
286--313,
DOI 10.4236/apm.2015.55029.
%
\bibitem{Ma29JA}
D. C. Mayer,
\textit{Index-\(p\) abelianization data of \(p\)-class tower groups},
29i\`emes Journ\'ees Arithm\'etiques,
Univ. of Debrecen (Debreceni Egyetem),
Hungary (Magyarorsz\'ag),
presentation July 09, 2015,
available from
\texttt{http://www.algebra.at/29JADebrecen.pdf}.
%
\bibitem{Ma2017}
D. C. Mayer,
\textit{Criteria for three-stage towers of \(p\)-class fields},
Adv. Pure Math.
\textbf{7}
(2015),
135--179,
DOI 10.4236/apm.2017.72008.
%
\bibitem{Ma2018a}
D. C. Mayer,
\textit{Deep transfers of \(p\)-class tower groups},
J. Appl. Math. Phys.
\textbf{6}
(2018),
36--50,
DOI 10.4236/jamp.2018.61005.
%
\bibitem{Ma2018b}
D. C. Mayer,
\textit{Annihilator ideals of two-generated metabelian \(p\)-groups},
J. Algebra Appl.
\textbf{17}
(2018),
no. 4,
DOI 10.1142/S0219498818500767.
%
\bibitem{Ne1989a}
B. Nebelung,
\textit{Klassifikation metabelscher \(3\)-Gruppen
mit Faktorkommutatorgruppe vom Typ \((3,3)\)
und Anwendung auf das Kapitulationsproblem},
Inauguraldissertation,
Universit\"at zu K\"oln,
1989.
%
\bibitem{Ne1989b}
B. Nebelung,
\textit{Anhang zu
Klassifikation metabelscher \(3\)-Gruppen
mit Faktorkommutatorgruppe vom Typ \((3,3)\)
und Anwendung auf das Kapitulationsproblem},
Band 2,
Universit\"at zu K\"oln,
1989.
%
\bibitem{Nm1977}
M. F. Newman,
\textit{Determination of groups of prime-power order},
pp. 73--84,
in: Group Theory, Canberra, 1975,
Lecture Notes in Math.,
Vol. \textbf{573},
Springer,
Berlin,
1977.
%
\bibitem{Ob1990}
E. A. O'Brien,
\textit{The \(p\)-group generation algorithm},
J. Symbolic Comput.
\textbf{9}
(1990),
677--698,
DOI 10.1016/S0747-7171(80)80082-X.
%
\bibitem{OEIS2025}
OEIS Foundation Inc. (N. J. A. Sloane),
\textit{The On-Line Encyclopedia of Integer Sequences} (OEIS),
2025,
Published electronically at 
\texttt{https://oeis.org/}.
%
\bibitem{So1933}
A. Scholz,
\textit{Idealklassen und Einheiten in kubischen K\"orpern},
Monatsh. Math. Phys.
\textbf{40}
(1933),
211--222.
%
\bibitem{SoTa1934}
A. Scholz und O. Taussky,
\textit{Die Hauptideale der kubischen Klassenk\"orper imagin\"ar quadratischer Zahlk\"orper:
ihre rechnerische Bestimmung und ihr Einflu\ss\ auf den Klassenk\"orperturm},
J. Reine Angew. Math.
\textbf{171}
(1934),
19--41.
%
\end{thebibliography}
\end{document}